\pdfoutput=1
\RequirePackage{ifpdf}
\ifpdf 
\documentclass[pdftex]{sigma}
\else
\documentclass{sigma}
\fi

\numberwithin{equation}{section}

\newtheorem{Theorem}{Theorem}[section]
\newtheorem*{Theorem*}{Theorem}
\newtheorem{Corollary}[Theorem]{Corollary}
\newtheorem{Lemma}[Theorem]{Lemma}
\newtheorem{Proposition}[Theorem]{Proposition}
 { \theoremstyle{definition}
\newtheorem{Definition}[Theorem]{Definition}

\newtheorem{Example}[Theorem]{Example}
\newtheorem{Remark}[Theorem]{Remark}
\newtheorem{Notation}[Theorem]{Notation}}

\input{xy}
\xyoption{all}

\newcommand{\eas}{\mathbf{EAS}}
\newcommand{\leas}{\ell\mathbf{EAS}}

\begin{document}

\allowdisplaybreaks

\newcommand{\arXivNumber}{2105.01326}

\renewcommand{\PaperNumber}{092}

\FirstPageHeading

\ShortArticleName{On Extended Associative Semigroups}

\ArticleName{On Extended Associative Semigroups}

\Author{Lo\"{\i}c FOISSY}

\AuthorNameForHeading{L.~Foissy}

\Address{Universit\'e Littoral C\^ote d'Opale, UR 2597 LMPA, Laboratoire de Math\'ematiques Pures\\ et Appliqu\'ees Joseph Liouville, 62100 Calais, France}
\Email{\href{mailto:foissy@univ-littoral.fr}{foissy@univ-littoral.fr}}
\URLaddress{\url{http://loic.foissy.free.fr/pageperso/accueil.html}}

\ArticleDates{Received June 12, 2025, in final form October 15, 2025; Published online October 26, 2025}

\Abstract{We study extended associative semigroups (briefly, EAS), an algebraic structure used to define generalizations of the operad of associative algebras, and the subclass of commutative extended diassociative semigroups (briefly, CEDS), which are used to define generalizations of the operad of pre-Lie algebras. We give families of examples based on semigroups or on groups, as well as a classification of EAS of cardinality two. We then define linear extended associative semigroups as linear maps satisfying a variation of the braid equation. We explore links between linear EAS and bialgebras and Hopf algebras. We also study the structure of non-degenerate finite CEDS and show that they are obtained by semi-direct and direct products involving two groups.}


\Classification{20M75; 16S10; 18M60; 16T05}

\section{Introduction}

It seems that the notion of \emph{family} parameterization of a given type of algebraic structure appeared firstly appears in the context of quantum field theory:
in \cite{Ebrahimi-Fard2007}, the authors introduced Rota--Baxter family algebras. This terminology is due to Li Guo \cite{Guo2009}. In the same spirit, family pre-Lie, dendriform or tridendriform algebras, among others have been introduced~\cite{Foissy49,Zhang2019,Zhang2020,Zhang2020-3,Manchon2020}. In all cases, the idea is to replace the operations defining the structure by a~bunch of operations indexed by a~semigroup $\Omega$;
the relations between the axioms are deformed using this structure on $\Omega$.
For example, if $(\Omega,\star)$ is a semigroup, an $(\Omega,\star)$-family associative algebra $A$ has a family $(*_\alpha)_{\alpha \in \Omega}$ of products, with the relations
\[x*_\alpha (y *_\beta z)=( x*_\alpha y)*_{\alpha \star \beta}z,\]
satisfied for any $\alpha,\beta,\gamma \in \Omega$ and any $x$, $y$, $z$ in $A$. In the same spirit, the notion of \emph{matching} parameterization can be used: for pre-Lie, it appears for example in the work of Bruned, Hairer and Zambotti
on regularity structures to solve stochastics PDEs \cite{Bruned2019,Bruned2023,Bruned2023-1,Foissy47}.
Matching Rota--Baxter algebras, associative, dendriform, pre-Lie algebras are introduced in \cite{Foissy47,Zhang2020}, see also \cite{Foissy49} for a two-parameter versions for pre-Lie algebras.
For example, a matching associative algebra has a family $(*_\alpha)_{\alpha \in \Omega}$ of products indexed by a set $\Omega$, with the relations
\[x*_\alpha (y *_\beta z)=(x *_\alpha y)*_\beta z.\]
Note that no specific structure is required on $\Omega$ in this case.
Attempts to unify these parameterizations have been done in \cite{Foissy56,Foissy48,Foissy50}. For example, for associative, following \cite{Foissy56},
given a set~$\Omega$ with two binary operations $\rightarrow$ and $\triangleright$, an $(\Omega,\rightarrow,\triangleright)$-associative algebra has a family $(*_\alpha)_{\alpha \in \Omega}$ of products, with the relations
\[x*_\alpha (y *_\beta z)=x *_{\alpha\triangleright \beta}(y*_{\alpha \rightarrow \beta}z).\]
Of course, usually these structures do not have any convenient property, and some conditions are imposed: roughly speaking, one imposes that the underlying combinatorics of the initial object is conserved, modulo a parameterization by $\Omega$.
This gives some constraints on $\Omega$. For associative algebras, $(\Omega,\rightarrow,\triangleright)$ has to be an extended associative semigroup (briefly, EAS):
\begin{gather*}
\forall \alpha,\beta,\gamma \in \Omega,\qquad \alpha \rightarrow (\beta \rightarrow \gamma)=(\alpha \rightarrow \beta) \rightarrow \gamma,\\
(\alpha \triangleright (\beta \rightarrow \gamma))\rightarrow (\beta \triangleright \gamma)=(\alpha \rightarrow \beta)\triangleright \gamma,\qquad
(\alpha \triangleright (\beta \rightarrow \gamma))\triangleright (\beta \triangleright \gamma)=\alpha \triangleright \beta.
\end{gather*}
In particular, $(\Omega,\rightarrow)$ is an associative semigroup. Here are some examples:
\begin{itemize}\itemsep=0pt
\item If $\Omega$ is a set, putting
$
 \forall \alpha,\beta \in \Omega$, $ \alpha \rightarrow\beta =\beta$, $\alpha \triangleright \beta =\alpha$,
we obtain an EAS, denoted by $\eas(\Omega)$. This EAS gives back matching associative algebras.
\item If $(\Omega,\rightarrow)$ is an associative semigroup, it is an EAS with
$\forall \alpha,\beta \in \Omega$, $\alpha \triangleright \beta =\alpha$.
This EAS is denoted by $\eas(\Omega,\rightarrow)$. It gives back $(\Omega,\rightarrow)$-family associative algebras.
\item If $(\Omega,\star)$ is a group, it is an EAS, with
$
 \forall \alpha,\beta \in \Omega$, $\alpha \rightarrow\beta =\beta$, $\alpha \triangleright \beta =\alpha\star \beta^{\star -1}$,
It is denoted by $\eas'(\Omega,\star)$. \end{itemize}
The two first examples explain why $\Omega$-matching and $(\Omega,\star)$-family associative algebras are very similar, in particular why the free objects are isomorphic as vector spaces:
this is fact works for the more general settings of $\Omega$-associative algebras over an EAS.
The same can be done with pre-Lie algebras, leading to the notion of commutative extended diassociative semigroup (briefly, CEDS). A CEDS is an EAS satisfying the complementary axioms
\begin{gather*}
\forall \alpha,\beta,\gamma \in \Omega,\qquad(\alpha \rightarrow \beta) \rightarrow \gamma =(\beta \rightarrow \alpha) \rightarrow \gamma,\qquad
\alpha \triangleright (\beta \rightarrow \gamma) =\alpha \triangleright \gamma,
\end{gather*}
The Koszul duality of quadratic operads applied to $\Omega$-pre-Lie algebras leads to the notion of dual CEDS, which are EAS with the complementary axioms
\begin{gather*}
 \forall \alpha,\beta,\gamma \in \Omega,\qquad (\alpha \triangleright \beta) \rightarrow \gamma=\alpha \rightarrow \gamma,\qquad
(\alpha \triangleright \beta) \triangleright \gamma =(\alpha \triangleright \gamma)\triangleright \beta.
\end{gather*}
For example, for any set $\Omega$, $\eas(\Omega)$ is both a CEDS and a dual CEDS. For any semigroup~$(\Omega,\star)$, $\eas(\Omega,\star)$ is a dual CEDS and is a CEDS if and only if
\begin{align*}
\forall \alpha,\beta,\gamma \in \Omega,\alpha \star (\beta \star \gamma)=(\alpha \star \beta) \star \gamma=(\beta \star \alpha) \star \gamma.
\end{align*}
For any group $(\Omega,\star)$, $\eas'(\Omega,\star)$ is a CEDS, and is a dual CEDS if and only if $\star$ is commutative.

The axioms of EAS can be reformulated using the maps
\begin{gather*}
\phi\colon\ \begin{cases}
\Omega^2\longrightarrow\Omega^2,\\
(\alpha,\beta)\longrightarrow(\alpha \rightarrow \beta,\alpha \triangleright \beta),
\end{cases}\qquad
\tau\colon\ \begin{cases}
\Omega^2\longrightarrow\Omega^2,\\
(\alpha,\beta)\longrightarrow(\beta,\alpha).
\end{cases}
\end{gather*}
Then $(\Omega,\rightarrow,\triangleright)$ is an EAS if and only if
\begin{align*}
(\operatorname{Id} \times \phi)\circ (\phi\times \operatorname{Id})\circ (\operatorname{Id} \times \phi)&=(\phi \times \operatorname{Id})\circ (\operatorname{Id} \times \tau)\circ (\phi\times \operatorname{Id}).
\end{align*}
Similar formulations can be done for CEDS and dual CEDS, see Lemma \ref{lem2.6}.
This reformulation naturally leads to the notion of linear EAS: an $\ell$EAS is pair $(A,\Phi)$, where $A$ is a vector space and $\Phi\colon A\otimes A\longrightarrow A\otimes A$ is a linear map satisfying the \emph{$\ell$EAS braid equation}
\begin{align*}
(\operatorname{Id} \otimes \Phi)\circ (\Phi\otimes \operatorname{Id})\circ (\operatorname{Id} \otimes \Phi)&=(\Phi \otimes \operatorname{Id})\circ (\operatorname{Id} \otimes \tau)\circ (\Phi\otimes \operatorname{Id}).
\end{align*}

In particular, let $(\Omega,\rightarrow,\triangleright)$ be a set with two operations. We denote by $\mathbb{K}\Omega$ the vector space generated by $\Omega$ and we define $\Phi\colon\mathbb{K}\Omega\otimes \mathbb{K}\Omega\longrightarrow\mathbb{K}\Omega\otimes \mathbb{K}\Omega$ by
\begin{align*}
&\forall \alpha,\beta\in \Omega,\qquad\Phi(\alpha\otimes \beta)=(\alpha \rightarrow \beta)\otimes (\alpha \triangleright \beta).
\end{align*}
Then $(\mathbb{K}\Omega,\Phi)$ is an $\ell$EAS if and only if $(\Omega,\rightarrow,\triangleright)$ is an EAS.
Not all the $\ell$EAS can be obtained in this way, see Example~\ref{ex4.1} for two-dimensional examples. Similar presentations of $\ell$CEDS and $\ell$CEDS can be established, see Definition \ref{defi4.1}.

The aim of this paper is a study of EAS, CEDS, and their linear versions. In the second section, after recalling the main definitions, we give a classification of EAS of cardinality 2,
which gives 13 non-isomorphic examples, 11 being CEDS, 7 being dual CEDS, 3 being non-degenerate (that is to say, with an invertible map $\phi$).
The third section is devoted to the study of non-degenerate finite CEDS. We prove that the three examples defined earlier are in fact fundamental bricks:
Theorem~\ref{theo3.14} states that any finite non-degenerate CEDS can be decomposed as the direct product of a semi-direct product $\eas(\Omega_1,*) \rtimes \eas'(\Omega_2,\star)$ with an~$\eas(\Omega_3)$,
where $(\Omega_1,*)$ is an abelian group, $(\Omega_2,\star)$ is a group and $\Omega_3$ is a nonempty set. The fourth section is devoted to linear versions of EAS.
We give firstly a family of 18 examples of~$\ell$EAS in dimension 2, then study the duality of $\ell$EAS (see Proposition \ref{prop4.2}), and left units, left counits and eigenvectors (see Definition \ref{defi4.3}).
If $(A,\Phi)$ is an $\ell$EAS, an element $a\in A$ is a~left unit if for any $b\in A$, $\Phi(a\otimes b)=b\otimes a$. An element $f \in A^*$ is a~left counit if for any $a,b\in A$, $(f\otimes \operatorname{Id})\circ \Phi(a\otimes b)=f(b)a$.
In particular, we characterize left units and counits and eigenvectors for linearization of non-degenerate CEDS in Proposition \ref{prop4.5}.
In the last section, we introduce two functors taking their values in the category of $\ell$EAS.
The first one (see Proposition \ref{prop5.1}) is defined on the category of bialgebras (not necessarily unitary nor counitary) and generalizes the construction of $\eas(\Omega,\rightarrow)$.
The second one (see Proposition \ref{prop5.7}) is defined on the category of Hopf algebras and generalizes the construction of $\eas'(\Omega,*)$.
In the case of an $\ell$EAS coming from a Hopf algebra, this is closely related to the notion of right integral (see Proposition \ref{prop5.11}).
We prove in Theorem~\ref{theo5.12} that we can associate to any convenient pair $(a,f)$ of a unit and a~counit a~bialgebra structure on $A$, recovering in this way $\ell$EAS coming from a bialgebra.
This is finally applied to $\ell$EAS defined from Hopf algebras of groups.

\begin{Notation}
$\mathbb{K}$ is a commutative field. All the vector spaces in this text are taken over $\mathbb{K}$.
\end{Notation}

\section{Extended (di)associative semigroups}

\subsection{Commutative extended diassociative semigroups}

Let us first recall this definition of \cite{Foissy49}, where it is related to a parameterization of the operad of dendriform algebras.

\begin{Definition}\qquad
\begin{enumerate}\itemsep=0pt
\item[(1)] A \emph{diassociative semigroup} is a family $(\Omega,\leftarrow,\rightarrow)$, where $\Omega$ is a nonempty set and $\leftarrow$, $\rightarrow\colon\Omega\times \Omega \longrightarrow \Omega$ are maps such that, for any $\alpha,\beta,\gamma\in \Omega$,
\begin{gather}
\label{eq1} (\alpha \leftarrow \beta)\leftarrow \gamma=\alpha\leftarrow (\beta \leftarrow \gamma)=\alpha \leftarrow(\beta \rightarrow \gamma),\\
 (\alpha \rightarrow \beta)\leftarrow \gamma=\alpha \rightarrow (\beta \leftarrow \gamma),\\
 (\alpha \rightarrow \beta)\rightarrow \gamma=(\alpha \leftarrow \beta) \rightarrow \gamma=\alpha \rightarrow (\beta \rightarrow \gamma).
\end{gather}
An \emph{extended diassociative semigroup} (briefly, EDS) is a family $(\Omega,\leftarrow,\rightarrow,\triangleleft,\triangleright)$, where $\Omega$ is a nonempty set
and $\leftarrow,\rightarrow,\triangleleft,\triangleright\colon\Omega\times \Omega \longrightarrow \Omega$ are maps such that
\begin{enumerate}\itemsep=0pt
\item $(\Omega,\leftarrow,\rightarrow)$ is a diassociative semigroup.
\item For any $\alpha,\beta,\gamma \in \Omega$,
\begin{gather}
\alpha\triangleright (\beta \leftarrow \gamma)=\alpha \triangleright \beta,\\
(\alpha \rightarrow \beta)\triangleleft \gamma=\beta \triangleleft \gamma, \\
(\alpha \triangleleft \beta)\leftarrow ((\alpha \leftarrow \beta) \triangleleft \gamma)=\alpha \triangleleft (\beta \leftarrow \gamma),\\
(\alpha \triangleleft \beta)\triangleleft ((\alpha \leftarrow \beta) \triangleleft \gamma)=\beta \triangleleft \gamma,
\\
(\alpha \triangleleft \beta) \rightarrow ((\alpha \leftarrow \beta) \triangleleft \gamma)=\alpha \triangleleft (\beta \rightarrow \gamma),\\
 (\alpha \triangleleft \beta) \triangleright ((\alpha \leftarrow \beta) \triangleleft \gamma)=\beta \triangleright \gamma,\\
(\alpha \triangleright (\beta \rightarrow \gamma))\leftarrow (\beta \triangleright \gamma)=(\alpha \leftarrow \beta) \triangleright \gamma,\\
(\alpha \triangleright (\beta \rightarrow \gamma))\triangleleft (\beta \triangleright \gamma)=\alpha \triangleleft \beta,\\
\label{eq12} (\alpha \triangleright (\beta \rightarrow \gamma))\rightarrow (\beta \triangleright \gamma)=(\alpha \rightarrow \beta)\triangleright \gamma,\\
\label{eq13} (\alpha \triangleright (\beta \rightarrow \gamma))\triangleright (\beta \triangleright \gamma)=\alpha \triangleright \beta.
\end{gather}
\end{enumerate}
\end{enumerate}
\end{Definition}

An EDS $(\Omega,\leftarrow,\rightarrow,\triangleleft,\triangleright)$ is commutative if for any $\alpha,\beta \in \Omega$,
\begin{align}
\label{eq14}\alpha \leftarrow \beta=\beta \rightarrow \alpha,\qquad\alpha \triangleleft \beta=\beta \triangleright \alpha.
\end{align}

Let us reformulate the definition of commutative EDS.

\begin{Proposition}
A commutative EDS $($briefly, CEDS$)$ is a triple $(\Omega,\rightarrow,\triangleright)$, where $\Omega$ is a~non\-empty set and $\rightarrow, \triangleright\colon\Omega^2\longrightarrow \Omega$ are maps such that, for any $\alpha,\beta,\gamma \in \Omega$,
\begin{gather}
\label{eq15} \alpha \rightarrow (\beta \rightarrow \gamma)=(\alpha \rightarrow \beta) \rightarrow \gamma=(\beta \rightarrow \alpha) \rightarrow \gamma,\\
\label{eq16} \alpha \triangleright (\beta \rightarrow \gamma)=\alpha \triangleright \gamma,\\
\label{eq17} (\alpha \triangleright \gamma)\rightarrow (\beta \triangleright \gamma)=(\alpha \rightarrow \beta)\triangleright \gamma,\\
\label{eq18} (\alpha \triangleright \gamma)\triangleright (\beta \triangleright \gamma)=\alpha \triangleright \beta.
\end{gather}\end{Proposition}

\begin{proof}
Replacing $\leftarrow$ and $\triangleleft$ in \eqref{eq1}--\eqref{eq13} with the help of \eqref{eq14}, we find \eqref{eq15}--\eqref{eq18}.
\end{proof}

\begin{Definition}[{\cite{Foissy56}}]
An \emph{extended associative semigroup} (briefly, EAS) is a triple $(\Omega,\rightarrow,\triangleright)$, where~$\Omega$ is a nonempty set and $\rightarrow,\triangleright\colon\Omega^2\longrightarrow \Omega$ are maps such that, for any $\alpha,\beta,\gamma \in \Omega$,
\begin{gather}
\label{eq19} \alpha \rightarrow (\beta \rightarrow \gamma) =(\alpha \rightarrow \beta) \rightarrow \gamma,\\
\tag{\ref{eq12}}(\alpha \triangleright (\beta \rightarrow \gamma))\rightarrow (\beta \triangleright \gamma) =(\alpha \rightarrow \beta)\triangleright \gamma,\\
\tag{\ref{eq13}}(\alpha \triangleright (\beta \rightarrow \gamma))\triangleright (\beta \triangleright \gamma) =\alpha \triangleright \beta.
\end{gather}\end{Definition}

\begin{Remark}
Let $(\Omega,\rightarrow,\leftarrow,\triangleright,\triangleleft)$ be an EDS. Then $(\Omega,\rightarrow,\triangleright)$ is an EAS, called the \emph{right part} of the EDS $(\Omega,\rightarrow,\leftarrow,\triangleright,\triangleleft)$.
We obtain a commutative triangle of functors
\[\xymatrix{\mathbf{CEDS}\ar@{^(->}[r]\ar@{^(->}[d]&\mathbf{EAS}.\\
\mathbf{EDS}\ar[ru]_{\text{\scriptsize{right part}}}&}\]
We shall see that not all the EAS are right parts of an EDS (see case \textbf{C6} in the classification of EAS of cardinality 2 in the next paragraph).
\end{Remark}

\begin{Example}\qquad
\begin{enumerate}\itemsep=0pt
\item[(1)] Let $\Omega$ be a set. We put
\begin{align*}
\forall (\alpha,\beta) \in \Omega^2,\qquad\begin{cases}
 \alpha \rightarrow \beta=\beta,\\
\alpha \triangleright \beta=\alpha.
\end{cases} \end{align*}
Then $(\Omega,\rightarrow,\triangleright)$ is an EAS, denoted by $\eas(\Omega)$. It is a CEDS.

\item[(2)] Let $(\Omega,\star)$ be an associative semigroup and let $\pi\colon\Omega\longrightarrow\Omega$ be an endomorphism of $(\Omega,\star)$ such that $\pi^2=\pi$. We put
$
\forall (\alpha,\beta) \in \Omega^2$, $ \alpha \triangleright \beta=\pi(\alpha)$.
It is an EAS, which we denote by~${\eas(\Omega,\star,\pi)}$. It is a CEDS if and only if for any~${\alpha,\beta,\gamma \in \Omega}$,
$(\alpha \star \beta)\star \gamma=(\beta \star \alpha)\star \gamma$.
We shall simply denote $\eas(\Omega,\star)$ instead of $\eas(\Omega,\star,\operatorname{Id}_\Omega)$. In particular, if $(\Omega,\star)$ is a~group, then $\eas(\Omega,\star)$ is a CEDS if and only if $(\Omega,\star)$ is abelian, which proves that not all EAS are CEDS.

\item[(3)] Let $\Omega$ be a set with an operation $\triangleright$ such that, for any $\alpha,\beta,\gamma\in \Omega$,
$(\alpha \triangleright \gamma)\triangleright (\beta\triangleright \gamma)=\alpha \triangleright \beta$.
We then put
$\forall (\alpha,\beta)\in \Omega^2$, $\alpha\rightarrow \beta =\beta$.
Then $(\Omega,\rightarrow,\triangleright)$ is a CEDS (so is an EAS). This holds, for example, if $(\Omega,\star)$ is an associative semigroup with the \emph{right inverse condition}
\begin{align*}
&\forall (\beta,\gamma)\in \Omega^2, \ \exists! \alpha \in \Omega,\qquad \alpha\star \beta=\gamma.
\end{align*}
This unique $\alpha$ is denoted by $\gamma \triangleright \beta$. Then, for any $\alpha,\beta,\gamma \in \Omega$,
\begin{align*}
((\alpha \triangleright \gamma)\triangleright (\beta\triangleright \gamma))\star \beta&=((\alpha \triangleright \gamma)\triangleright (\beta\triangleright \gamma))\star ((\beta \triangleright \gamma)\star \gamma)\\
&=(((\alpha \triangleright \gamma)\triangleright (\beta\triangleright \gamma))\star (\beta \triangleright \gamma))\star \gamma
=(\alpha \triangleright \gamma)\star \gamma=\alpha,
\end{align*}
so $(\alpha \triangleright \gamma)\triangleright (\beta\triangleright \gamma)=\alpha \triangleright \beta$. This EAS is denoted by $\eas'(\Omega,\star)$.
The right inverse condition holds for example if $(\Omega,\star)$ is a group, and then
$\alpha \triangleright \beta=\alpha \star \beta^{\star-1}$.
It also holds for semigroups which are not groups. For example, if $\Omega$ is a nonempty set, we give it an associative product defined by
$\forall \alpha,\beta\in \Omega$, $ \alpha \star \beta=\alpha$.
It satisfies the right inverse condition and, for any $\alpha,\beta \in \Omega$, $\alpha \triangleright \beta=\alpha$. Note that for this example, $\eas'(\Omega,\star)=\eas(\Omega)$.
\end{enumerate}\end{Example}

\begin{Definition} \label{defi2.4}
Let $(\Omega,\rightarrow,\triangleright)$ be an EAS. We shall say that it is \emph{non-degenerate} if the following map is bijective:
\[\phi\colon \begin{cases}
\Omega^2\longrightarrow\Omega^2,\\
(\alpha,\beta)\longrightarrow(\alpha\rightarrow\beta,\alpha \triangleright\beta).
\end{cases}\]
\end{Definition}

If $\Omega$ is a non-degenerate EAS, the structure implied on $\Omega$ by $\phi^{-1}$ will be studied in the next paragraph.

\begin{Example}\qquad
\begin{enumerate}\itemsep=0pt
\item[(1)] Let $\Omega$ be a set. In $\eas(\Omega)$, for any $\alpha,\beta \in \Omega$, $\phi(\alpha,\beta)=(\beta,\alpha)$, so $\eas(\Omega)$ is non-degenerate.
\item[(2)] Let $(\Omega,\star)$ be a group. Then $\eas(\Omega,\star)$ is non-degenerate. Indeed, in this case, $\phi(\alpha,\beta)=(\alpha\star\beta,\alpha)$, so $\phi$ is a bijection, of inverse given by $\phi^{-1}(\alpha,\beta)=\bigl(\beta,\beta^{\star-1}\star \alpha\bigr)$.
\item[(3)] Let $(\Omega,\star)$ be an associative semigroup with the right inverse condition. Then $\eas'(\Omega,\star)$ is non-degenerate.
Indeed, in this case, $\phi(\alpha,\beta)=(\beta,\alpha\triangleright \beta)$, so $\phi$ is a bijection, of inverse given by $\phi^{-1}(\alpha,\beta)=(\beta\star \alpha,\alpha)$.
\end{enumerate}\end{Example}

\subsection{Dual commutative extended semigroups}

\begin{Definition}
Let $(\Omega,\rightarrow,\triangleright)$ be a set with two binary operations. We shall say that it is a~\emph{dual CEDS} if, for any $\alpha,\beta,\gamma \in \Omega$,
\begin{gather}
\tag{\ref{eq19}} (\alpha \rightarrow \beta)\rightarrow \gamma=\alpha \rightarrow (\beta \rightarrow \gamma),\\
\tag{\ref{eq12}}(\alpha \triangleright (\beta \rightarrow \gamma))\rightarrow (\beta \triangleright \gamma)=(\alpha \rightarrow \beta) \triangleright \gamma,\\
\tag{\ref{eq13}} (\alpha \triangleright (\beta \rightarrow \gamma))\triangleright (\beta \triangleright \gamma)=\alpha \triangleright \beta,\\
 (\alpha \triangleright \beta) \rightarrow \gamma=\alpha \rightarrow \gamma,\\
 (\alpha \triangleright \beta) \triangleright \gamma=(\alpha \triangleright \gamma)\triangleright \beta.
\end{gather}\end{Definition}

\begin{Example}\qquad
\begin{enumerate}\itemsep=0pt
 \item[(1)] If $\Omega$ is a set, then $\eas(\Omega)$ is a dual CEDS.
 \item[(2)] If $(\Omega,\star)$ is a semigroup and $\pi\colon\Omega\longrightarrow \Omega$ is a semigroup morphism such that $\pi^2=\pi$, then~${\eas(\Omega,\star,\pi)}$ is a dual CEDS if and only if
$
 \forall \alpha,\beta \in \Omega$, $\pi(\alpha)\star \beta =\alpha \star \beta$.
 In particular, $\eas(\Omega,\star)$ is a dual CEDS.
 \item[(3)] If $(\Omega,\star)$ is a semigroup with the right inverse condition, then $\eas'(\Omega,\star)$ is a dual CEDS if and only if
$\forall \alpha,\beta,\gamma \in \Omega$, $(\alpha \triangleright \beta) \triangleright \gamma =(\alpha \triangleright \gamma) \triangleright \beta$.
 This is equivalent to
$\forall \alpha,\beta,\gamma \in \Omega$, $\alpha \star \beta \star \gamma =\alpha \star \gamma \star \beta$.
 In the case where $(\Omega,\star)$ is a group, $\eas'(\Omega,\star)$ is a dual CEDS if and only if $(\Omega,\star)$ is abelian.
 \end{enumerate}\end{Example}

The following lemma, proved in \cite{Foissy48}, is a reformulation of the axioms of EAS, CEDS and dual CEDS with the help of the map $\phi$.

\begin{Lemma}\label{lem2.6}
Let $(\Omega,\rightarrow,\triangleright)$ be a set with two binary operations. We consider the maps
\begin{align*}
\phi\colon\ \begin{cases}
\Omega^2 \longrightarrow \Omega^2,\\
(\alpha,\beta) \longrightarrow (\alpha \rightarrow \beta,\alpha \triangleright \beta),
\end{cases} \qquad
\tau\colon\ \begin{cases}
\Omega^2 \longrightarrow \Omega^2,\\
(\alpha,\beta) \longrightarrow (\beta,\alpha).
\end{cases}
\end{align*}
Then
\begin{enumerate}\itemsep=0pt
\item[$(1)$] $(\Omega,\rightarrow,\triangleright)$ is an EAS if and only if
\begin{align}
\label{eq23}(\operatorname{Id} \times \phi)\circ (\phi\times \operatorname{Id})\circ (\operatorname{Id} \times \phi) =(\phi \times \operatorname{Id})\circ (\operatorname{Id} \times \tau)\circ (\phi\times \operatorname{Id}).
\end{align}
\item[$(2)$] $(\Omega,\rightarrow,\triangleright)$ is a CEDS if and only if
\begin{gather}
(\operatorname{Id} \times \phi)\circ (\phi\times \operatorname{Id})\circ (\operatorname{Id} \times \phi) =(\phi \times \operatorname{Id})\circ (\operatorname{Id} \times \tau)\circ (\phi\times \operatorname{Id}),\tag{\ref{eq23}}\\
(\operatorname{Id} \times \phi)\circ (\operatorname{Id} \times \tau) \circ (\tau \times \operatorname{Id})\circ (\phi \times \operatorname{Id})\nonumber\\
\qquad =(\tau \times \operatorname{Id})\circ (\phi \times \operatorname{Id})\circ (\operatorname{Id} \times \phi)\circ (\operatorname{Id} \times \tau).\label{eq24}
\end{gather}
$(\Omega,\rightarrow,\triangleright)$ is a dual CEDS if and only if
\begin{gather}
(\operatorname{Id} \times \phi)\circ (\phi\times \operatorname{Id})\circ (\operatorname{Id} \times \phi) =(\phi \times \operatorname{Id})\circ (\operatorname{Id} \times \tau)\circ (\phi\times \operatorname{Id}),\tag{\ref{eq23}}\\
(\phi \times \operatorname{Id})\circ (\tau \times \operatorname{Id})\circ (\operatorname{Id} \times \tau)\circ (\operatorname{Id} \times \phi) \nonumber\\
\qquad=(\operatorname{Id} \times \tau)\circ (\operatorname{Id} \times \phi)\circ (\phi \times \operatorname{Id})\circ (\tau \times \operatorname{Id}).\label{eq25}
\end{gather}\end{enumerate}\end{Lemma}

With this reformulation, the following result becomes immediate, as the
inversion of \eqref{eq23} gives \eqref{eq23} again and the inversion of \eqref{eq24} gives \eqref{eq25}.

\begin{Proposition}\label{prop2.7}
Let $(\Omega, \rightarrow,\triangleright)$ be a set with two binary operations. We shall say that $(\Omega,\rightarrow, \triangleright)$ is non-degenerate if the map $\phi$ of Definition~{\rm\ref{defi2.4}} is a bijection. If so, we put
\[\phi^{-1}\colon \ \begin{cases}
\Omega^2 \longrightarrow \Omega^2,\\
(\alpha,\beta) \longrightarrow (\alpha\curvearrowright \beta,\alpha \blacktriangleright \beta).
\end{cases}\]
Then $(\Omega,\rightarrow,\triangleright)$ is an EAS $($resp.\ a CEDS, a dual CEDS$)$ if and only if $(\Omega,\curvearrowright,\blacktriangleright)$ is an EAS $($resp.\ a dual CEDS, a CEDS$)$.
\end{Proposition}

\subsection{EAS of cardinality two}

Here is a classification of EAS of cardinality two, which we obtained by an exhaustive study of the $2^8$ possibilities of pairs of operations. The underlying set is $\Omega=\{X,Y\}$ and the products will be given by the pair of matrices
\begin{align*}
\begin{pmatrix}
X\rightarrow X&X\rightarrow Y\\Y\rightarrow X&Y\rightarrow Y
\end{pmatrix},\qquad
\begin{pmatrix}
X\triangleright X&X\triangleright Y\\Y\triangleright X&Y\triangleright Y
\end{pmatrix}.
\end{align*}
We shall use the two maps
\begin{align*}
\pi_X\colon \ \begin{cases}
\Omega\longrightarrow\Omega,\\
\alpha\longrightarrow X,
\end{cases}\qquad
\pi_Y\colon \ \begin{cases}
\Omega\longrightarrow\Omega,\\
\alpha\longrightarrow Y.
\end{cases}
\end{align*}
We respect the indexation of EDS of \cite{Foissy49}.
$$\renewcommand{\arraystretch}{1.1}
\begin{array}{|c|c|c|l|l|}
\hline \text{Case}&\rightarrow&\triangleright&\text{Description}&\text{Comments}\\
\hline \hline \mathbf{A1}&\begin{pmatrix}
X&X\\X&X
\end{pmatrix}&
\begin{pmatrix}
X&X\\X&X
\end{pmatrix}&\eas(\Omega,\rightarrow,\pi_X)&\begin{array}{@{}l@{}}\text{CEDS, dual CEDS,}\\
 \text{right part of \textbf{D1}}\end{array}\\
\hline \mathbf{A2}&\begin{pmatrix}
X&X\\X&X
\end{pmatrix}&
\begin{pmatrix}
X&X\\Y&Y
\end{pmatrix}&\eas(\Omega,\rightarrow)&\begin{array}{@{}l@{}}\text{CEDS, dual CEDS,}\\
\text{right part of \textbf{D2}}\end{array}\\
\hline \mathbf{C1} &\begin{pmatrix}
X&X\\X&Y
\end{pmatrix}&
\begin{pmatrix}
X&X\\X&X
\end{pmatrix}&\eas(\Omega,\rightarrow,\pi_X)&\text{CEDS, right part of \textbf{C4}}\\
\hline \mathbf{C3}&\begin{pmatrix}
X&X\\X&Y
\end{pmatrix}&
\begin{pmatrix}
X&X\\Y&Y
\end{pmatrix}&\eas(\mathbb{Z}/2\mathbb{Z},\times)&\text{CEDS, dual CEDS}\\
\hline \mathbf{C5} &\begin{pmatrix}
X&X\\X&Y
\end{pmatrix}&
\begin{pmatrix}
Y&Y\\Y&Y
\end{pmatrix}&\eas((\mathbb{Z}/2\mathbb{Z},\times),\pi_Y)&\text{CEDS, right part of \textbf{C2}}\\
\hline \mathbf{C6}&\begin{pmatrix}
X&X\\X&Y
\end{pmatrix}&
\begin{pmatrix}
X&X\\Y&X
\end{pmatrix}&&\\
\hline \mathbf{E1'}\text{--}\mathbf{E2'}&\begin{pmatrix}
X&X\\Y&Y
\end{pmatrix}&
\begin{pmatrix}
X&X\\X&X
\end{pmatrix}&\eas(\Omega,\rightarrow,\pi_X)&\text{right part of \textbf{E1} and \textbf{E2}}\\
\hline \mathbf{E3'}&\begin{pmatrix}
X&X\\Y&Y
\end{pmatrix}&
\begin{pmatrix}
X&X\\Y&Y
\end{pmatrix}&\eas(\Omega,\rightarrow)&\begin{array}{@{}l@{}}\text{dual CEDS},\\
\text{right part of \textbf{E3}}\end{array}\\
\hline\mathbf{ F1}&\begin{pmatrix}
X&Y\\X&Y
\end{pmatrix}&
\begin{pmatrix}
X&X\\X&X
\end{pmatrix}&\eas(\Omega,\rightarrow,\pi_X)&\begin{array}{@{}l@{}}\text{CEDS, dual CEDS},\\
\text{right part of \textbf{B1}, \textbf{F2}, \textbf{G1}}\\
\text{and \textbf{G2}}\end{array}\\
\hline \mathbf{F3}&\begin{pmatrix}
X&Y\\X&Y
\end{pmatrix}&
\begin{pmatrix}
X&X\\Y&Y
\end{pmatrix}&\eas(\Omega)&\begin{array}{@{}l@{}}\text{CEDS, dual CEDS},\\
\text{non-degenerate},\\
\text{right part of \textbf{B2} and \textbf{G3}}\end{array}\\
\hline \mathbf{F4}&\begin{pmatrix}
X&Y\\X&Y
\end{pmatrix}&
\begin{pmatrix}
X&Y\\Y&X
\end{pmatrix}&\eas'(\mathbb{Z}/2\mathbb{Z},+)&\begin{array}{@{}l@{}}\text{CEDS, dual CEDS},\\
\text{non-degenerate},\\
\text{right part of \textbf{F5}}\end{array}\\
\hline \mathbf{H1}&\begin{pmatrix}
X&Y\\Y&X
\end{pmatrix}&
\begin{pmatrix}
X&X\\X&X
\end{pmatrix}&\eas(\mathbb{Z}/2\mathbb{Z},+,\pi_X)&\text{CEDS}\\
\hline \mathbf{H2}&\begin{pmatrix}
X&Y\\Y&X
\end{pmatrix}&
\begin{pmatrix}
X&X\\Y&Y
\end{pmatrix}&\eas(\mathbb{Z}/2\mathbb{Z},+)&\begin{array}{@{}l@{}}\text{CEDS, dual CEDS},\\
\text{non-degenerate}\end{array}\\
\hline \end{array}
$$
For the cases \textbf{C3}, \textbf{C5}, \textbf{F4}, \textbf{H1} and \textbf{H2}, $\Omega$ is identified with $\mathbb{Z}/2\mathbb{Z}$, $X$ being $\overline{0}$ and $Y$ being $\overline{1}$.

\begin{Remark}
With similar methods, it is possible to prove that there are three non-degenerate EAS of cardinality 3 up to isomorphism: $\eas(\{1,2,3\})$, $\eas(\mathbb{Z}/3\mathbb{Z},+)$ and $\eas'(\mathbb{Z}/3\mathbb{Z},+)$. All of them are both CEDS and dual CEDS.
\end{Remark}

\section{Structure of non-degenerate finite CEDS}

We now turn to CEDS, and prove the structure Theorem~\ref{theo3.14} after several intermediate results.

\subsection{Preliminary results}

\begin{Lemma}\label{lem3.1}
let $\Omega$ be a finite non-degenerate EAS.
\begin{enumerate}\itemsep=0pt
\item[$(1)$] Let $\Omega'$ be a sub-EAS of $\Omega$. Then $\Omega'$ is non-degenerate.
\item[$(2)$] Let $\sim$ be an equivalence on $\Omega$, compatible with the EAS structure. Then the quotient EAS $\Omega/{\sim}$ is non-degenerate.
\end{enumerate}
\end{Lemma}

\begin{proof}
(1) By restriction, $\phi_{\Omega'}=(\phi_\Omega)_{\mid \Omega'^2}$ is injective. As $\Omega'$ is finite, it is a bijection. So $\Omega'$ is non degenerate.

(2) Let $\pi\colon\Omega\longrightarrow \Omega/{\sim}$ be the canonical surjection. Then $\phi_{\Omega/{\sim}}\circ \pi=(\pi\otimes \pi)\circ \phi_\Omega$. As $\phi$ is surjective, $\phi_{\Omega/{\sim}}$ is surjective.
As $\Omega/{\sim}$ is finite, it is a bijection. So $\Omega/{\sim}$ is non-degenerate.
 \end{proof}

\begin{Definition}
Let $(\Omega,\rightarrow,\triangleright)$ be an EAS. For any $\alpha \in \Omega$, we put
\begin{align*}
\phi_\alpha\colon\ \begin{cases}
\Omega \longrightarrow \Omega,\\
\beta \longrightarrow \alpha\rightarrow \beta,
\end{cases}\qquad
\psi_\alpha\colon\ \begin{cases}
\Omega \longrightarrow \Omega,\\
\beta \longrightarrow \beta\triangleright\alpha.
\end{cases}
\end{align*}
We shall say that $(\Omega,\rightarrow,\triangleright)$ is strongly non-degenerate if for any $\alpha\in \Omega$, $\phi_\alpha$ is bijective.
\end{Definition}

\begin{Remark}
As the product $\rightarrow$ is associative, for any $\alpha,\beta \in \Omega$,
$\phi_\alpha \circ \phi_\beta=\phi_{\alpha \rightarrow \beta}$.
\end{Remark}

\begin{Lemma} \label{lem3.3}
Let $(\Omega,*)$ be an associative semigroup. The following conditions are equivalent:
\begin{enumerate}\itemsep=0pt
\item[$(1)$] $\eas(\Omega,*)$ is non-degenerate.
\item[$(2)$] $\eas(\Omega,*^{\rm op})$ is strongly non-degenerate.
\item[$(3)$] $(\Omega,*^{\rm op})$ has the right inverse condition.
\end{enumerate}
\end{Lemma}

\begin{proof}
Let $\alpha$, $\beta$, $\gamma$, $\delta\in \Omega$. Then
\begin{align*}
\phi(\alpha,\beta)=(\gamma,\delta) \Longleftrightarrow \begin{cases}
\alpha*\beta=\gamma,\\
\alpha=\delta.
\end{cases}\end{align*}
So
\begin{align*}
\phi\text{ is bijective}&\Longleftrightarrow\forall (\gamma,\delta)\in \Omega^2,
\exists! \beta \in \Omega, \delta*\beta=\gamma\\
&\Longleftrightarrow \text{in $\eas(\Omega,*^{\rm op})$, }\forall \delta \in \Omega, \phi_{\delta}\text{ is bijective}\\
&\Longleftrightarrow\text{$(\Omega,*^{\rm op})$ has the right inverse condition}. \tag*{\qed}
\end{align*} \renewcommand{\qed}{}\end{proof}

\begin{Lemma}\label{lem3.4}
Let $(\Omega,\rightarrow,\triangleright)$ be a finite non-degenerate CEDS. Then it is strongly non-degenerate.
\end{Lemma}

\begin{proof}
Let $\alpha,\gamma,\gamma' \in \Omega$ such that $\phi_\alpha(\gamma)=\phi_\alpha(\gamma')$. In other words, $\alpha \rightarrow \gamma=\alpha\rightarrow \gamma'$. By \eqref{eq16},
\begin{align*}
\alpha \triangleright\gamma&=\alpha \triangleright (\alpha \rightarrow \gamma)=\alpha \triangleright (\alpha \rightarrow \gamma')=\alpha \triangleright \gamma'.
\end{align*}
Therefore, $\phi(\alpha,\gamma)=\phi(\alpha,\gamma')$. As $\phi$ is injective, $\gamma=\gamma'$, so $\phi_\alpha$ is injective. As $\Omega$ is finite, $\phi_\alpha$ is bijective.
\end{proof}

\begin{Lemma}\label{lem3.5}
Let $\Omega=(\Omega,\rightarrow,\triangleright)$ be a non-degenerate EAS, such that
$\forall \alpha,\beta \in \Omega$, $\alpha \rightarrow \beta =\beta$.
There exists a product $*$ on $\Omega$, making it a semigroup with the right inverse condition, such that~${\Omega=\eas'(\Omega,*)}$. For any $\beta\in \Omega$, $\psi_\beta$ is bijective and its inverse is
\begin{align*}
\phi'_\beta\colon\ \begin{cases}
\Omega \longrightarrow \Omega,\\
\alpha \longrightarrow \alpha*\beta.
\end{cases}
\end{align*}
Moreover, for any $\beta,\gamma \in \Omega$,
\begin{align}
\label{eq26}
\psi_\beta \circ \psi_\gamma =\psi_{\beta*\gamma},\qquad \psi_{\beta \triangleright \gamma} =\psi_\beta \circ \psi_\gamma^{-1}.
\end{align}\end{Lemma}

\begin{proof}
Note that for any $\alpha \in \Omega$, $\phi_\alpha=\operatorname{Id}_\Omega$. Let $\alpha,\beta,\gamma,\delta\in \Omega$. Then
\begin{align*}
\phi(\alpha,\beta)=(\gamma,\delta) \Longleftrightarrow \begin{cases}
\beta=\gamma,\\
\alpha\triangleright \beta=\delta.
\end{cases}\end{align*}
Hence,
\begin{align*}
\phi\text{ is bijective}&\Longleftrightarrow\forall (\gamma,\delta)\in \Omega^2, \exists! \alpha \in \Omega, \alpha \triangleright\gamma=\delta\Longleftrightarrow\forall\gamma \in \Omega, \psi_\gamma\text{ is bijective}.
\end{align*}
Putting $\phi^{-1}(\alpha,\beta)=(\alpha\curvearrowright \beta,\alpha \blacktriangleright \beta)$, by Proposition \ref{prop2.7} $(\Omega,\curvearrowright,\blacktriangleright)$ is an EAS, so $\curvearrowright$ is associative.
Moreover, $\phi^{-1}(\alpha,\beta)=(\alpha \curvearrowright \beta,\alpha)$, so $(\Omega,\curvearrowright,\blacktriangleleft)=\eas(\Omega,\curvearrowright)$.
By Lemma \ref{lem3.3}, if $*=\curvearrowright^{\rm op}$, then $*$ has the right inverse condition. Moreover, for any $\alpha,\beta \in \Omega$,
\begin{align*}
\phi^{-1}\circ \phi(\alpha,\beta)&=\phi^{-1}(\beta,\alpha \triangleright \beta)=((\alpha \triangleright \beta)*\beta,\beta)=(\alpha,\beta).
\end{align*}
Hence, the unique element $\gamma \in \Omega$ such that $\gamma*\beta=\alpha$ is $\alpha \triangleright \beta$: consequently, $\Omega=\eas'(\Omega,*)$. Moreover, for any $\alpha,\beta \in \Omega$,
$\phi'_\beta \circ \psi_\beta(\alpha)=(\alpha\triangleright \beta) *\beta=\alpha$.
So $\phi'_\beta\circ \psi_\beta=\operatorname{Id}_\Omega$. As $\psi_\beta$ is bijective, $\psi_\beta^{-1}=\phi'_\beta$.

Let $\beta,\gamma\in \Omega$. Then, for any $\alpha \in \Omega$,
$\phi'_\gamma \circ \phi'_\beta(\alpha)=\alpha *\beta*\gamma=\phi'_{\beta*\gamma}$.
So $\phi'_\gamma\circ \phi'_\beta=\phi'_{\beta*\gamma}$. Inverting, $\psi_\beta\circ \psi_\gamma=\psi_{\beta*\gamma}$. As a consequence,
$\psi_{\beta \triangleright \gamma}\circ \psi_\gamma=\psi_{(\beta \triangleright \gamma)*\gamma}=\psi_\beta$,
which induces the last formula. \end{proof}

\begin{Lemma}\label{lem3.6}
Let $\Omega=(\Omega,\rightarrow,\triangleright)$ be a non-degenerate EAS such that $\forall \alpha,\beta \in \Omega$,
$\alpha \rightarrow\beta=\beta$.
Then $\Omega_\psi=\{\psi_\alpha, \alpha\in \Omega\}$ is a subgroup of the group of permutations of $\Omega$.
\end{Lemma}

\begin{proof} Direct consequence of \eqref{eq26}. \end{proof}

\begin{Proposition}\label{prop3.7}
Let $\Omega=\eas'(\Omega,*)$, where $(\Omega,*)$ is a finite semigroup with the right inverse condition. We define an equivalence $\sim$ on $\Omega$ by $\alpha \sim \beta$ if $\psi_\alpha=\psi_\beta$. Then
\begin{enumerate}\itemsep=0pt
\item[$(1)$] $\sim$ is compatible with the EAS structure of $\Omega$. Therefore, $\Omega/{\sim}$ is an EAS.
\item[$(2)$] There exists a product $\star$ on $\Omega/{\sim}$, making it a group, such that $\Omega/{\sim}=\eas'(\Omega/{\sim},\star)$.
\item[$(3)$] There exists a sub-EAS $\Omega_0$ of $\Omega$, such that the restriction to $\Omega_0$ of the canonical surjection~${\pi\colon\Omega\longrightarrow \Omega/{\sim}}$ is an isomorphism.
\end{enumerate}\end{Proposition}

\begin{proof}
(1) Let $\alpha,\beta \in \Omega$, such that $\alpha \sim\beta$. Then $\psi_\alpha=\psi_\beta$. Let $\gamma \in \Omega$. Then $\alpha \rightarrow \gamma=\beta \rightarrow \gamma=\gamma$, and $\gamma \rightarrow \alpha=\alpha \sim \beta=\gamma \rightarrow \beta$.
As $\psi_\alpha=\psi_\beta$, $\gamma \triangleright \alpha=\gamma \triangleright \beta$. Moreover, by Lemma \ref{lem3.5},
\[\psi_{\alpha \triangleright \gamma}=\psi_\alpha \circ \psi_\gamma^{-1}=\psi_\beta\circ \psi_\gamma^{-1}=\psi_{\beta \triangleright \gamma},\]
so $\alpha \triangleright \gamma\sim \beta \triangleright \gamma$: $\sim$ is compatible with the EAS structure.

(2) By Lemma \ref{lem3.1}, $\Omega/{\sim}$ is non-degenerate. By Lemma \ref{lem3.5}, there exists a product $\star$ satisfying the right inverse condition, such that $\Omega/{\sim}=\eas'(\Omega/{\sim},\star)$. We consider the map
\[\psi\colon\ \begin{cases}
(\Omega/{\sim},\star) \longrightarrow (\mathfrak{S}_{\Omega/{\sim}},\circ),\\
\overline{\alpha} \longrightarrow \psi_{\overline{\alpha}}.
\end{cases}\]
By Lemma \ref{lem3.5}, this is a semigroup morphism. Let us prove that it is injective. We assume that $\psi_{\overline{\alpha}}=\psi_{\overline{\beta}}$.
In other words, for any $\gamma \in \Omega$, $\gamma \triangleright \alpha\sim \gamma \triangleright \beta$, or equivalently, $\psi_{\gamma \triangleright \alpha}=\psi_{\gamma \triangleright \beta}$. Moreover,
\[\psi_{\gamma \triangleright \alpha}=\psi_\gamma \circ \psi_\alpha^{-1}=\psi_{\gamma \triangleright \beta}=\psi_\gamma \circ \psi_\beta^{-1}.\]
As $\psi_\gamma$ is bijective, $\psi_\alpha=\psi_\beta$, so $\overline{\alpha}=\overline{\beta}$.

By Lemma \ref{lem3.6}, there exists $e\in \Omega/{\sim}$, such that $\psi_e=\operatorname{Id}_{\Omega/{\sim}}$. For any $\overline{\alpha}\in \Omega/{\sim}$,
\[\psi_{e\triangleright \overline{\alpha}}=\psi_e\circ \psi_{\overline{\alpha}}^{-1}=\psi_{\overline{\alpha}}^{-1},\]
so $\psi(\Omega/{\sim})$ is a subgroup of $\mathfrak{S}_{\Omega/{\sim}}$. Consequently, $(\Omega/{\sim},\star)$ is a group.

(3) By Lemma \ref{lem3.6}, there exists $\beta_0\in \Omega$ such that $\psi_{\beta_0}=\operatorname{Id}_\Omega$. We put
\[\Omega_0=\{\beta_0\triangleright \alpha,\, \alpha \in \Omega\}=\{\psi_\alpha(\beta_0),\, \alpha \in \Omega\}.\]
As the product $\rightarrow$ of $\Omega$ is trivial, this is a sub-semigroup of $(\Omega,\rightarrow)$. Let $\beta_0\triangleright \alpha$, $\beta_0\triangleright \beta \in \Omega_0$.
\[(\beta_0\triangleright \alpha)\triangleright (\beta_0\triangleright \beta)=\psi_{\beta_0\triangleright \gamma}\circ \psi_{\alpha}(\beta_0)=\psi_{(\beta_0\triangleright \gamma)*\alpha}(\beta_0)\in \Omega_0,\]
so $\Omega_0$ is a sub-EAS of $\Omega$.

Let us assume that $\beta_0\triangleright \alpha\sim \beta_0\triangleright \beta$. Then
\begin{align*}
\psi_{\beta_0\triangleright \alpha}=\psi_{\beta_0}\circ \psi_\alpha^{-1}=\psi_{\alpha}^{-1}
=\psi_{\beta_0\triangleright \beta}=\psi_{\beta_0}\circ \psi_\beta^{-1}=\psi_{\beta}^{-1},
\end{align*}
so $\psi_\alpha=\psi_\beta$. Hence, $\beta_0\triangleright\alpha=\beta_0\triangleright \beta$, which proves that \smash{$\pi_{\mid \Omega_0}$} is injective.
By Lemma \ref{lem3.6}, there exists $\beta \in \Omega$ such that $\psi_\beta=\psi_\alpha^{-1}$. We consider $\beta_0\triangleright \beta \in \Omega_0$. Then
$\psi_{\beta_0\triangleright \beta}=\psi_{\beta_0}\circ \psi_\beta^{-1}=\psi_\alpha$,
so $\beta_0\triangleright \beta\sim \alpha$. Hence, $\pi_{\mid \Omega_0}$ is surjective.
\end{proof}

\begin{Theorem} \label{theo3.8}
Let $\Omega=\eas'(\Omega,*)$, where $(\Omega,*)$ is a finite semigroup with the right inverse condition. There exists a group $(\Omega_0,\star)$ and a set $\Omega_1$ such that
$\Omega\approx \eas(\Omega_1)\times \eas'(\Omega_0,\star)$.
\end{Theorem}

\begin{proof}
We keep the notations of the proof of Proposition \ref{prop3.7}. As the sub-EAS $\Omega_0$ is isomorphic to $\Omega/{\sim}$, it is a group for the law $*$, and $\Omega_0=\eas'(\Omega_0,*)$.
Let $e$ be the unit of the group~$(\Omega/{\sim},\star)$. We consider
$\Omega_1=\{\alpha\in \Omega,\, \overline{\alpha}=e\}$.
Let us prove that $\Omega_1=\{\alpha \in \Omega,\, \psi_\alpha=\operatorname{Id}_\Omega\}$.

$\supseteq$: if $\psi_\alpha=\operatorname{Id}_\Omega$, for any $\beta \in \Omega$, $\overline{\beta} \star \overline{\alpha}^{*-1}=\overline{\psi_\alpha(\beta)}=\overline{\beta}$, so $\overline{\alpha}=e$ and $\alpha \in \Omega_1$.

$\subseteq$: if $\overline{\alpha}=e$, then for any $\beta \in \Omega$, $\overline{\beta \triangleright \alpha}=\overline{\beta}$, so $\beta \triangleright \alpha\sim \beta$: in other words, $\psi_{\beta \triangleright \alpha}=\psi_\beta$. Then
\[\psi_{\beta \triangleright \alpha}=\psi_\beta \circ \psi_{\alpha}^{-1}=\psi_\beta.\]
As $\psi_\beta$ is bijective, $\psi_\alpha=\operatorname{Id}_\Omega$.

Therefore, for any $\alpha \in \Omega$, for any $\beta \in \Omega_1$,
$\alpha \triangleright \beta=\psi_\beta(\alpha)=\alpha$.
As a consequence, $\Omega_1=\eas(\Omega_1)$. We consider the map
\[\theta\colon\ \begin{cases}
\Omega_1\times \Omega_0 \longrightarrow \Omega,\\
(\alpha,\beta) \longrightarrow \alpha*\beta.
\end{cases}\]
Let us prove that $\theta$ is injective. If $\theta(\alpha,\beta)=\theta(\alpha',\beta')$, in $\Omega/{\sim}$,
$\overline{\alpha}\star \overline{\beta}=\overline{\beta}=\overline{\alpha'}\star \overline{\beta'}=\overline{\beta'}$.
As $\pi_{\mid \Omega_0}$ is injective, $\beta=\beta'$. Because of the right inverse condition for $*$, $\alpha=\alpha'$.

Let us prove that $\theta$ is surjective. Let $\gamma \in \Omega$. There exists a unique $\beta\in \Omega_0$ such that $\psi_{\overline{\gamma}}=\psi_{\overline{\beta}}$.
We put $\alpha=\gamma \triangleright \beta$, so $\gamma=\alpha*\beta$. Moreover,
\smash{$\psi_{\overline{\alpha}}=\psi_{\overline{\gamma}}\circ \psi_{\overline{\beta}}^{-1}=\operatorname{Id}_{\Omega/{\sim}}$}.
So \smash{$e\triangleright \overline{\alpha}=\overline{\alpha}^{\star -1}=\psi_{\overline{\alpha}}(e)=e$}, and finally $\alpha\in \Omega_1$.

Let $(\alpha,\beta)$ and $(\alpha',\beta')\in \Omega_1 \times \Omega_0$. In $\Omega$, as $\alpha'\in \Omega_1$,
$\alpha*\beta*\alpha'*\beta'=\alpha*(\beta*\beta')$,
which implies that~${(\alpha*\beta)\triangleright (\alpha'*\beta')=\alpha* \bigl(\beta*\beta'^{*-1}\bigr)}$.
So $\theta$ is an isomorphism of EAS from $\eas(\Omega_1)\times \eas'(\Omega_0,*)$ to $\Omega$.
\end{proof}

\subsection{Non-degenerate finite CEDS}

\begin{Lemma}\label{lem3.9}
Let $(\Omega,\rightarrow,\triangleright)$ be a strongly non-degenerate finite EAS. Then $\Omega_\phi=(\{\phi_\alpha, \alpha\in \Omega\},\circ)$ is a group. The following map is a surjective morphism of semigroups:
\[\phi\colon \ \begin{cases}
(\Omega,\rightarrow) \longrightarrow \Omega_\phi,\\
\alpha \longrightarrow \phi_\alpha.
\end{cases}\]
\end{Lemma}

\begin{proof}
We already observed that for any $\alpha,\beta \in \Omega$, $\phi_\alpha\circ \phi_\beta=\phi_{\alpha\rightarrow \beta}$, so $\phi$ is a semigroup morphism.
By hypothesis, for any $\alpha \in \Omega$, $\phi_\alpha$ is a bijection, so belongs to the symmetric group~$\mathfrak{S}_\Omega$ of permutations of $\Omega$.
As $\Omega$ is finite, for any $\alpha \in \Omega$, there exists $n\geqslant 2$ such that~${\phi_\alpha^n=\operatorname{Id}_\Omega}$.
Then $\phi_{\alpha^{\rightarrow n}}=\operatorname{Id}_\Omega$, so $\Omega_\phi$ is a monoid. Putting $\beta=\alpha^{\rightarrow(n-1)}$, $\phi_\beta\circ \phi_\alpha=\phi_\alpha\circ \phi_\beta=\operatorname{Id}_\Omega$, so $\Omega_\phi$ is a group.
\end{proof}

\begin{Proposition}\label{prop3.10}
Let $\Omega=(\Omega,\rightarrow,\triangleright)$ be a finite non-degenerate EAS, such that for any $\alpha \in \Omega$, $\phi_\alpha$ is a bijection. We put
$
\Omega^\rightarrow=\{\alpha\in \Omega, \phi_\alpha=\operatorname{Id}_\Omega\}$, $\Omega^\triangleright =\{\beta \in \Omega, \psi_\beta=\operatorname{Id}_\Omega\}$.
Then
\begin{enumerate}\itemsep=0pt
\item[$(1)$] $\Omega^\rightarrow$ is a non-degenerate sub-EAS of $\Omega$.
\item[$(2)$] If $\Omega^\triangleright$ is nonempty, it is a non-degenerate sub-EAS of $\Omega$.
\item[$(3)$] If $\Omega^\triangleright$ is nonempty, then $\Omega^\triangleright \cap \Omega^\rightarrow$ is nonempty.
\item[$(4)$] If $\Omega$ is a CEDS, $\Omega^\triangleright$ is nonempty.
\end{enumerate}\end{Proposition}

\begin{proof}
(1) Recall that for any $\alpha,\beta \in \Omega$, $\phi_\alpha\circ \phi_\beta=\phi_{\alpha \rightarrow \beta}$. This easily implies that $\Omega^\rightarrow$ is stable under $\rightarrow$.
By Lemma \ref{lem3.9}, there exists $\alpha \in \Omega$, such that $\phi_\alpha=\operatorname{Id}_\Omega$, so $\Omega^\rightarrow$ is nonempty.

Let $\alpha,\beta \in \Omega^\rightarrow$. Let us consider $\gamma \in \Omega$. As $\phi$ is bijective, there exist $\beta',\gamma'\in \Omega$ such that
\[(\beta'\rightarrow \gamma',\beta'\triangleright \gamma')=(\beta,\gamma).\]
Then
\begin{align*}
\phi_{\alpha \triangleright \beta}(\gamma)&=(\alpha\triangleright \beta)\rightarrow \gamma=\bigl(\alpha \triangleright\bigl(\beta'\rightarrow \gamma'\bigr)\bigr)\rightarrow \bigl(\beta' \triangleright \gamma'\bigr)=\beta' \triangleright \gamma'=\gamma.
\end{align*}
So $\phi_{\alpha \triangleright \beta}=\operatorname{Id}_\Omega$ and $\alpha \triangleright \beta \in \Omega^\rightarrow$. By Lemma \ref{lem3.1}, $\Omega^\rightarrow$ is a non-degenerate sub-EAS.

(2) Let $\beta,\gamma \in \Omega^\triangleright$. As $\psi_\gamma=\operatorname{Id}_\Omega$, $\beta\triangleright \gamma=\beta \in \Omega^\triangleright$. For any $\alpha \in \Omega$, by \eqref{eq12} and \eqref{eq13},
\begin{gather*}
(\alpha \triangleright (\beta\rightarrow \gamma))\rightarrow (\beta \triangleright \gamma) =(\alpha \rightarrow \beta)\triangleright \gamma=\alpha \rightarrow \beta,\\
(\alpha \triangleright (\beta\rightarrow \gamma))\triangleright (\beta \triangleright \gamma) =\alpha \triangleright \beta.
\end{gather*}
So $\phi(\alpha \triangleright (\beta\rightarrow \gamma),\beta \triangleright \gamma)=\phi(\alpha,\beta)$.
As $\phi$ is injective, $\alpha \triangleright (\beta\rightarrow \gamma)=\alpha$, so $\psi_{\beta \rightarrow \gamma}=\operatorname{Id}_\Omega$ and~${\beta \rightarrow \gamma\in \Omega^\triangleright}$.
If $\Omega^\triangleright$ is nonempty, by Lemma \ref{lem3.1}, it is non-degenerate.

(3) Let us take $\alpha \in \Omega^\triangleright$. The permutation $\phi_\alpha$ is of finite order as $\Omega$ is finite, so there exists~${n\geqslant 2}$, such that $\phi_{\alpha}^n=\phi_{\alpha^{\rightarrow n}}=\operatorname{Id}_\Omega$.
Putting $\beta=\alpha^{\rightarrow n}$, then $\beta \in \Omega^{\triangleright}$ (as it is a sub-CEDS) and $\beta \in \Omega^\rightarrow$ as $\phi_\beta=\operatorname{Id}_\Omega$.

(4) Let us consider the EAS associated to the inverse of $\phi$ (see Proposition \ref{prop2.7}), which we denote by $(\Omega,\curvearrowright,\blacktriangleright)$.
By the first point, there exists $\alpha \in \Omega$ such that for any $\beta \in \Omega$, $\alpha \curvearrowright \beta=\beta$. In other words, for any $\beta \in \Omega$,
$\phi^{-1}(\alpha,\beta)=(\beta,\alpha \blacktriangleright \beta)$.
This implies that $\phi_\beta(\alpha \blacktriangleright \beta)=\alpha$. The inverse of the bijection $\phi_\beta$ is the map
\[\phi'_\beta\colon \ \begin{cases}
\Omega \longrightarrow \Omega,\\
\alpha \longrightarrow \alpha \blacktriangleright \beta.
\end{cases}\]
As $\Omega$ is finite, there exists $\beta'\in \Omega$ such that $\phi_\beta^{-1}=\phi_{\beta'}$. Hence,
\[\beta=\beta \triangleright (\alpha \blacktriangleright \beta)=\beta\triangleright \bigl(\beta'\rightarrow \alpha\bigr)=\beta \triangleright \alpha,\]
by \eqref{eq16}. So $\alpha \in \Omega^\triangleright$. \end{proof}

\begin{Proposition}\label{prop3.11}
Let $(\Omega,\rightarrow,\triangleright)$ be a finite non-degenerate CEDS.
\begin{enumerate}\itemsep=0pt
\item[$(1)$] We define an equivalence $\equiv$ on $\Omega$ by
$\beta \equiv \beta'$ if $\exists \alpha \in \Omega$, $\beta' =\alpha\rightarrow \beta$.
This equivalence is compatible with the CEDS structure. Therefore, $\Omega/\equiv$ is a CEDS.
\item[$(2)$] The restriction to $\Omega^\rightarrow$ of the canonical surjection $\pi\colon\Omega\longrightarrow \Omega/\equiv$ is an isomorphism.
\item[$(3)$] $\Omega=\Omega^\triangleright\rightarrow \Omega^\rightarrow$.
\end{enumerate}\end{Proposition}

\begin{proof}
(1) The relation $\equiv$ can be reformulated as: there exists $\phi_\alpha\in \Omega_\phi$, such that $\phi_\alpha(\beta)=\beta'$. By Lemmas \ref{lem3.4} and \ref{lem3.9}, $\Omega_\phi$ is a group.
This easily implies that $\equiv$ is an equivalence: its classes are the orbits of the action of the group $\Omega_\phi$ over $\Omega$.

Let us assume that $\beta \equiv \beta'$: we put $\beta'=\alpha\rightarrow \beta$. Let $\gamma \in \Omega$. Then $\gamma \rightarrow \beta\equiv \beta\equiv \beta'\equiv \gamma \rightarrow \beta'$ by definition of $\equiv$.
Moreover, $\beta' \rightarrow \gamma=\alpha \rightarrow (\beta \rightarrow \gamma)$, so $\beta'\rightarrow \gamma\equiv \beta \rightarrow \gamma$. By \eqref{eq13},
\begin{gather*}
\beta'\triangleright \gamma=(\alpha \rightarrow \gamma)\triangleright \gamma=(\alpha \triangleright (\beta \rightarrow \gamma))\rightarrow (\beta \triangleright \gamma)\equiv \beta \triangleright \gamma,\\
\gamma \triangleright \beta'=\gamma \triangleright(\alpha \rightarrow \beta)=\gamma \triangleright \beta.
\end{gather*}
So $\equiv$ is compatible with the CEDS structure.

(2) Let $\alpha \in \Omega$. As $\phi_\alpha$ is bijective, there exists a unique $\beta \in \Omega$ such that $\alpha \rightarrow \beta=\alpha$. Then
$\phi_\alpha=\phi_{\alpha \rightarrow \beta}=\phi_\alpha\circ \phi_\beta$.
As $\phi_\alpha$ is bijective, $\phi_\beta=\operatorname{Id}_\Omega$, so $\beta \in \Omega^\rightarrow$ and $\alpha \equiv \beta$. This proves that~$\pi_{\mid \Omega^\rightarrow}$ is surjective.

Let $\beta,\beta'\in \Omega^\rightarrow$, such that $\beta\equiv \beta'$. There exists $\alpha\in \Omega$ such that $\alpha\rightarrow\beta=\beta'$. Then
\[\operatorname{Id}_\Omega=\phi_{\beta'}=\phi_{\alpha}\circ \phi_{\beta}=\phi_\alpha,\]
so $\phi_\alpha=\operatorname{Id}_\Omega$. We deduce that $\beta'=\phi_\alpha(\beta)=\beta$. Hence, $\pi_{\mid \Omega^\rightarrow}$ is injective.

(3) By Proposition \ref{prop3.10}, there exists $\beta_0\in \Omega^\triangleright \cap \Omega^\rightarrow$.
Let $\beta \in \Omega$. As $\pi_{\mid \Omega^\rightarrow}$ is bijective, there exists $\beta_1\in \Omega^\rightarrow$, such that $\beta\equiv \beta_1$. We put $\beta=\alpha \rightarrow \beta_1$. As $\beta_0\in \Omega^\rightarrow$,
$\beta=\alpha \rightarrow {\beta_0\rightarrow \beta_1}$. Moreover, for any $\gamma \in \Omega$, as $\beta_0\in \Omega^\triangleright$, by \eqref{eq16},
$\gamma \triangleright (\alpha \rightarrow \beta_0)=\gamma \triangleright \beta_0=\gamma$,
so $\alpha \rightarrow \beta_0\in \Omega^\triangleright$. \end{proof}

\begin{Proposition}\label{prop3.12}
Let $(\Omega,\rightarrow,\triangleright)$ be a finite non-degenerate CEDS. We define an equivalence on $\Omega^\triangleright$ by
\[\alpha' \sim \alpha'' \Longleftrightarrow \exists \alpha \in \Omega^\rightarrow, \ \alpha''=\alpha'\rightarrow \alpha.\]
\begin{enumerate}\itemsep=0pt
\item[$(1)$] This equivalence is compatible with the CEDS structure, and $\Omega'=\Omega^\triangleright/{\sim}$ is a non-degen\-er\-ate CEDS. Moreover, $(\Omega',\rightarrow)$ is an abelian group and $\Omega'=\eas(\Omega',\rightarrow)$.
\item[$(2)$] The following map is a semigroup isomorphism:
\[\theta\colon \ \begin{cases}
(\Omega' \times \Omega^\rightarrow,\rightarrow) \longrightarrow (\Omega,\rightarrow),\\
(\overline{\alpha},\beta) \longrightarrow \alpha \rightarrow \beta.
\end{cases}\] \end{enumerate} \end{Proposition}

\begin{proof} We firstly introduce an auxiliary map, defined by
\[\Theta\colon \ \begin{cases}
\Omega^\triangleright \times \Omega^\rightarrow \longrightarrow \Omega,\\
(\alpha,\beta) \longrightarrow \alpha \rightarrow \beta.
\end{cases}\]
By Proposition \ref{prop3.11}, it is surjective. Let us prove that $\Theta(\alpha',\beta')=\Theta(\alpha'',\beta'')$ if and only if $\alpha'\sim \alpha''$ and $\beta'=\beta''$.

Let us assume that $\Theta(\alpha',\beta')=\Theta(\alpha'',\beta'')$. As $\phi_{\alpha'}$ is bijective, there exists $\alpha \in \Omega$, $\alpha''=\alpha'\rightarrow \alpha$.
As $\alpha'\rightarrow \beta'=\alpha''\rightarrow \beta''$ and $\beta',\beta''\in \Omega^\rightarrow$,
\[\phi_{\alpha'}=\phi_{\alpha'}\circ \phi_{\beta'}=\phi_{\alpha'\rightarrow \beta'}=\phi_{\alpha''\rightarrow \beta''}=\phi_{\alpha''}.\]
Hence,
$\phi_{\alpha''}=\phi_{\alpha'}\circ \phi_{\alpha}=\phi_{\alpha'}$.
As $\phi_{\alpha'}$ is bijective, $\phi_\alpha=\operatorname{Id}_\Omega$, so $\alpha \in \Omega^\rightarrow$: we obtain that $\alpha'\equiv \alpha''$. As $\phi_{\alpha'}=\phi_{\alpha''}$, $\phi_{\alpha'}(\beta')=\phi_{\alpha''}(\beta'')=\phi_\alpha'(\beta'')$.
As $\phi_{\alpha'}$ is injective, $\beta=\beta'$. Conversely, if $\alpha \in \Omega^\rightarrow$, $\alpha'\rightarrow \alpha \rightarrow\beta'=\alpha'\rightarrow \beta'$.

As a consequence, $\sim$ is indeed an equivalence, $\theta$ is well-defined and is a bijection. It remains to show that $\sim$ is compatible with the CEDS structure of $\Omega^\triangleright$.
Let $\alpha',\alpha''\in \Omega^\triangleright$, such that $\alpha' \sim \alpha''$. We put $\alpha''=\alpha'\rightarrow \alpha$, with $\alpha \in \Omega^\rightarrow$. Let $\beta \in \Omega^\triangleright$. Then
\begin{align*}
\alpha'\triangleright \beta=\alpha'\sim \alpha''=\alpha''\triangleright \beta,\qquad\beta \triangleright \alpha'=\beta=\beta \triangleright \alpha''.
\end{align*}
Moreover,
\begin{align*}
\alpha''\rightarrow \beta=\alpha''\rightarrow \alpha \rightarrow \beta=\alpha''\rightarrow \beta,\qquad\beta \rightarrow \alpha''=\beta \rightarrow \alpha'\rightarrow \alpha\sim \beta \rightarrow \alpha'.
\end{align*}
Therefore, $\Omega^\triangleright/{\sim}$ is a CEDS. By Lemma \ref{lem3.1}, it is non-degenerate.

Let $\alpha,\alpha'\in \Omega^\triangleright$ and $\beta,\beta'\in \Omega^\rightarrow$. As $\beta\in \Omega^\rightarrow$,
\begin{align*}
\theta(\overline{\alpha},\beta)\rightarrow \theta\bigl(\overline{\alpha'},\beta'\bigr)&=\alpha \rightarrow \beta \rightarrow \alpha' \rightarrow \beta'=\alpha \rightarrow \alpha' \rightarrow \beta'\\
&=\alpha \rightarrow \alpha' \rightarrow \beta \rightarrow\beta'=\theta\bigl(\overline{\alpha}\rightarrow \overline{\alpha'},\beta \rightarrow \beta'\bigr).
\end{align*}
So $\theta$ is an isomorphism for the products $\rightarrow$.

Let us now study the CEDS $\Omega'$. By definition of $\Omega^\triangleright$, for any $\overline{\alpha}$, $\overline{\beta}\in \Omega'$, $\overline{\alpha}\triangleright \overline{\beta}=\overline{\alpha}$, so $\Omega'=\eas'(\Omega',\rightarrow)$.
By Proposition \ref{prop3.10}, $\Omega^\rightarrow$ is nonempty. Let us prove that
\[\Omega'^\rightarrow=\big\{\overline{\alpha}, \alpha\in \Omega^\triangleright\cap \Omega^\rightarrow\big\}.\]
$\supseteq$ is obvious. Let us take $\overline{\alpha}\in \Omega'^\rightarrow$. Then, for any $\beta \in \Omega^\triangleright$, $\overline{\alpha}\rightarrow\overline{\beta}=\overline{\beta}$:
there exists $\gamma \in \Omega^\rightarrow$, $\alpha \rightarrow \beta=\beta \rightarrow \gamma$. Therefore,
$\phi_\alpha\circ \phi_\beta=\phi_\beta \circ \phi_\gamma=\phi_\beta$,
as $\phi_\gamma=\operatorname{Id}_\Omega$. As $\phi_\beta$ is a bijection, $\phi_\alpha=\operatorname{Id}_\Omega$, so $\alpha \in \Omega^\rightarrow$. Let $\alpha,\beta \in \Omega^\triangleright\cap \Omega^\rightarrow$.
As $\phi_\alpha$ is bijective, there exists $\beta'\in \Omega$, $\alpha \rightarrow \beta'=\beta$. Then
\[ \operatorname{Id}_\Omega=\phi_\beta=\phi_\alpha\circ \phi_{\beta'}=\phi_{\beta'},\]
so $\beta'\in \Omega^\rightarrow$: we proved that $\alpha \sim \beta$. As a conclusion, there exists a unique $\overline{e}\in \Omega'$, such that for any $\overline{\alpha}\in \Omega'$, $\overline{e}\rightarrow \overline{\alpha}=\overline{\alpha}$.

Let us choose $\overline{\alpha}\in \Omega'$. As $\phi_\alpha$ is bijective, there exists $\overline{e'}\in \Omega'$ such that $\overline{\alpha}\rightarrow \overline{e'}=\overline{\alpha}$.
Let~${\overline{\beta}\in \Omega'}$. Then $\overline{\alpha}\rightarrow \overline{e'}\rightarrow \overline{\beta}=\overline{\alpha}\rightarrow \overline{\beta}$:
in other words, $\alpha \rightarrow e'\rightarrow \beta\sim \alpha \rightarrow \beta$, and there exists~${\gamma \in \Omega^\rightarrow}$, such that $\alpha \rightarrow e'\rightarrow \beta=\alpha \rightarrow \beta \rightarrow \gamma$.
As $\phi_\alpha$ is injective, $e'\rightarrow \beta=\beta \rightarrow \gamma\sim \beta$, so~${\overline{e'}\rightarrow \overline{\beta}=\overline{\beta}}$ for any $\overline{\beta}\in \Omega'$.
By unicity of $\overline{e}$, $\overline{e'}=\overline{e}$: for any $\overline{\alpha}\in \Omega'$, $\overline{\alpha}\rightarrow\overline{e}=\overline{\alpha}$, so $\overline{e}$ is a unit of $(\Omega',\rightarrow)$.
By~\eqref{eq15}, for $\gamma=\overline{e}$, we deduce that $(\Omega',\rightarrow)$ is an abelian monoid.
Let $\overline{\alpha}\in \Omega'$. As~${\phi_\alpha}$ is surjective, there exists $\overline{\alpha}'\in \Omega'$ such that $\overline{\alpha}\rightarrow \overline{\alpha}'=\overline{e}$. So $(\Omega',\rightarrow)$ is a group. \end{proof}

{\samepage
\begin{Proposition}
Let $(\Omega,*)$ be an associative semigroup such that for any $\alpha,\beta,\gamma \in \Omega$,
\[\alpha*\beta * \gamma=\beta *\alpha*\gamma.\]
Let $(\Omega',\rightarrow,\triangleright)$ be a CEDS, and $\prec\colon\Omega \times \Omega'\longrightarrow \Omega$ be a map such that for any $\alpha,\beta \in \Omega$, for any~${\beta',\gamma'\in \Omega'}$,
\begin{gather}
\label{eq27} \alpha \prec \bigl(\beta'\rightarrow \gamma'\bigr)=\alpha \prec \gamma',\\
\label{eq28} \bigl(\alpha *\beta\bigr)\prec \gamma' =\bigl(\alpha \prec \gamma'\bigr)*\bigl(\beta \prec \gamma'\bigr),\\
\label{eq29} \bigl(\alpha\prec \gamma'\bigr)\prec \bigl(\beta'\triangleright \gamma'\bigr)=\alpha \prec \beta',
\end{gather}
we define two products $\rightarrow$ and $\triangleright$ on $\Omega\times \Omega'$ in the following way: for any $(\alpha,\alpha')$, $(\beta,\beta')\in \Omega \times \Omega'$,
\begin{gather*}
\bigl(\alpha,\alpha'\bigr)\rightarrow \bigl(\beta,\beta'\bigr)=\bigl(\alpha*\beta,\alpha'\rightarrow \beta'\bigr),\qquad
\bigl(\alpha,\alpha'\bigr)\triangleright \bigl(\beta,\beta'\bigr)=\bigl(\alpha \prec \beta',\alpha'\triangleright \beta'\bigr).
\end{gather*}
Then $(\Omega\times \Omega,'\rightarrow,\triangleright)$ is a CEDS, which we denote by $\Omega \rtimes_\prec \Omega'$.
\end{Proposition}

\begin{proof}
Direct verifications.
\end{proof}

}

\begin{Remark}
If for any $(\alpha,\alpha')\in \Omega\times \Omega'$, $\alpha \prec \alpha'=\alpha$, we recover the direct product $\Omega \times \Omega'$ of EAS.
\end{Remark}

\begin{Theorem}\label{theo3.14}
Let $\Omega$ be a finite non-degenerate CEDS. There exist an abelian group $(\Omega_1,*)$, a group $(\Omega_2,\star)$, a left action $\succ\colon\Omega_2\times \Omega_1\longrightarrow \Omega_1$ of $(\Omega_2,\star)$ on $(\Omega_1,*)$
by group automorphisms, and a nonempty set $\Omega_3$ such that $\Omega$ is of the form
\[(\eas(\Omega_1,*) \rtimes_\succ \eas'(\Omega_2,\star))\times \eas(\Omega_3),\]
with the products given by
\begin{gather*}
(\alpha_1,\alpha_2,\alpha_3)\rightarrow(\beta_1,\beta_2,\beta_3)=(\alpha_1*\beta_1,\beta_2,\beta_3),\\
(\alpha_1,\alpha_2,\alpha_3)\triangleright(\beta_1,\beta_2,\beta_3)=\bigl(\beta_2\succ \alpha_1,\alpha_2\star \beta_2^{\star-1},\alpha_3\bigr).
\end{gather*}
\end{Theorem}

\begin{proof}
Let us consider the map $\theta$ of Proposition \ref{prop3.12}. For any $\alpha,\alpha'\in \Omega^\triangleright$, and for any $\beta,\beta'\in \Omega^\rightarrow$, by \eqref{eq16} and \eqref{eq17},
\begin{align*}
(\alpha \rightarrow \beta) \triangleright \bigl(\alpha'\rightarrow \beta'\bigr)&=(\alpha \rightarrow \beta)\triangleright \beta'=\bigl(\alpha \triangleright \beta'\bigr)\rightarrow \bigl(\beta \triangleright \beta'\bigr).
\end{align*}
Let $\beta_0\in \Omega^\triangleright \cap \Omega^\rightarrow$. Then, as $\beta_0\in \Omega^\rightarrow$,
\[(\alpha \rightarrow \beta) \triangleright \bigl(\alpha'\rightarrow \beta'\bigr)=\underbrace{\bigl(\alpha \triangleright \beta'\bigr)\rightarrow \beta_0}_{=\gamma_1}\rightarrow\underbrace{\bigl(\beta \triangleright \beta'\bigr)}_{=\gamma_2}.\]
Obviously, $\gamma_2\in \Omega^\rightarrow$. For any $\gamma \in \Omega$, by \eqref{eq16},
$\gamma \triangleright \gamma_1=\gamma \triangleright \beta_0=\gamma$,
so $\gamma_1\in \Omega^\triangleright$. We then put, for any $\overline{\alpha}\in \Omega'$, and for any $\beta\in \Omega^\rightarrow$,
$\overline{\alpha}\prec \beta=\overline{\alpha\triangleright \beta}\rightarrow \overline{\beta_0}$.
Then, for any $\overline{\alpha}, \overline{\alpha'}\in \Omega'$, and for any~${\beta,\gamma \in \Omega^\rightarrow}$,
$\theta(\overline{\alpha},\beta)\triangleright \theta\bigl(\overline{\alpha'},\beta'\bigr)=\theta\bigl(\overline{\alpha}\prec \beta',\beta \triangleright \beta'\bigr)$.
Then
\begin{align*}
\overline{\alpha}\prec(\beta \rightarrow \gamma)&=\overline{\alpha \triangleright(\beta \rightarrow \gamma)}\rightarrow \overline{\beta_0}=\overline{\alpha \triangleright \gamma}\rightarrow \overline{\beta_0}=\overline{\alpha}\prec \gamma,
\end{align*}
which proves \eqref{eq27}. As $\beta_0\in \Omega^\rightarrow$,
\begin{align*}
(\overline{\alpha}\rightarrow \overline{\beta})\prec \gamma&=\overline{(\alpha \rightarrow \beta) \triangleright \gamma}\rightarrow \overline{\beta_0}=\overline{\alpha \triangleright \gamma}\rightarrow\overline{\beta \triangleright \gamma}\rightarrow \overline{\beta_0}\\
&=\overline{\alpha \triangleright \gamma}\rightarrow \overline{\beta_0}\rightarrow\overline{\beta \triangleright \gamma}\rightarrow \overline{\beta_0}=(\overline{\alpha}\prec \gamma)\rightarrow (\overline{\beta}\prec \gamma),
\end{align*}
which proves \eqref{eq28}. Then
\begin{align*}
(\overline{\alpha}\prec \gamma)\prec (\beta \triangleright\gamma)&=\overline{(\alpha \triangleright \gamma)\triangleright (\beta \triangleright \gamma)}\rightarrow \overline{\beta_0\triangleright(\beta \triangleright \gamma)}\rightarrow \overline{\beta_0}=\overline{\alpha \triangleright \beta}\rightarrow \overline{\beta_0\triangleright(\beta \triangleright \gamma)}\rightarrow \overline{\beta_0}\\
&=\overline{\alpha \triangleright \beta}\rightarrow \overline{\beta_0}=\overline{\alpha}\prec \beta,
\end{align*}
which proves \eqref{eq29}. For the last equality, we used that $\beta_0\triangleright(\beta \triangleright \gamma)\in \Omega^\rightarrow$, as $\beta$, $\beta_0$ and $\gamma$ belong to $\Omega^\rightarrow$.

We finally obtain that $\theta$ is an isomorphism between $\Omega'\rtimes_\prec \Omega^\rightarrow$ and $\Omega$. We put $\Omega'=\eas(\Omega_1,*)$.
From Theorem~\ref{theo3.8}, we obtain a decomposition of $\Omega^\rightarrow$ of the form $\eas'(\Omega_2,\star)\times \eas(\Omega_3)$.
The map $\prec\colon\Omega_1\times \Omega_2\times \Omega_3\longrightarrow \Omega_1$ satisfies \eqref{eq27}--\eqref{eq29}. In this particular case, \eqref{eq27} becomes trivial, and \eqref{eq28}, \eqref{eq29} can be reformulated in this way:
for any $\alpha_1,\beta_1\in \Omega_1$, $\beta_2,\gamma_2\in \Omega_2$, $\beta_3,\gamma_3\in \Omega_3$,
\begin{gather*}
(\alpha_1*\beta_1)\prec(\gamma_2,\gamma_3)=(\alpha_1\prec (\gamma_2,\gamma_3))*(\beta_1\prec (\gamma_2,\gamma_3)),\\
(\alpha_1\prec (\gamma_2,\gamma_3))\prec(\beta_2,\beta_3)=\alpha_1\prec (\beta_2\star \gamma_2,\beta_3).
\end{gather*}
The products of $\Omega$ are given in this way: for any $\alpha_i,\beta_i\in \Omega_i$, with $1\leqslant i\leqslant 3$,
\begin{gather*}
\begin{split}
& (\alpha_1,\alpha_2,\alpha_3)\rightarrow(\beta_1,\beta_2,\beta_3)=(\alpha_1*\beta_1,\beta_2,\beta_3),\\
& (\alpha_1,\alpha_2,\alpha_3)\triangleright(\beta_1,\beta_2,\beta_3)=\bigl(\alpha_1\prec (\beta_2,\beta_3),\alpha_2\star \beta_2^{\star-1},\alpha_3\bigr).
\end{split}
\end{gather*}
For any $(\beta_2,\beta_3)\in \Omega_2\times \Omega_3$, we consider
\[\psi^\prec_{\beta_2,\beta_3}\colon \ \begin{cases}
\Omega_1 \longrightarrow \Omega_1,\\
\alpha_1 \longrightarrow \alpha_1\prec (\beta_2,\beta_3).
\end{cases}\]
As $\Omega$ is non-degenerate, necessarily $\psi_{\beta_2,\beta_3}$ is injective. As $\Omega$ is finite, $\psi_{\beta_2,\beta_3}$ is a bijection.
Moreover, by \eqref{eq30}, for any $(\beta_2,\beta_3)$, $(\gamma_2,\gamma_3)\in \Omega_2\times \Omega_3$,
$\psi^\prec_{\beta_2,\beta_3}\circ \psi^\prec_{\gamma_2,\gamma_3}=\psi^\prec_{\beta_2\star \gamma_2,\beta_3}$.
For $\beta_2=\gamma_2$ being the unit $e$ of $\Omega_2$ and $\beta_3=\gamma_3$, we obtain that \smash{$\bigl(\psi^\prec_{e,\beta_3}\bigr)^2=\psi^\prec_{e,\beta_3}$}.
As it is a~bijection, $\psi^\prec_{e,\beta_3}=\operatorname{Id}_{\Omega_1}$ for any $\beta_3\in \Omega_3$. Hence,
\[\psi^\prec_{e,\beta_3}\circ \psi^\prec_{\gamma_2,\gamma_3}=\psi^\prec_{\gamma_2,\gamma_3}=\psi^\prec_{\gamma_2,\beta_3},\]
so $\psi_{\beta_2,\beta_3}$ does not depend on $\beta_3$. We denote this map by $\psi_{\beta_2}$. Note that we proved that $\psi_{\beta_{e_{\Omega_2}}}=\operatorname{Id}_{\Omega_1}$.
We put, for any $\alpha_1\in \Omega_1$, $\beta_2\in \Omega_2$, $\alpha_1\prec \beta_2=\psi^\prec_{\beta_2}(\alpha_1)$. We finally obtain that the products in $\Omega$ are given by
\begin{gather*}
(\alpha_1,\alpha_2,\alpha_3)\rightarrow(\beta_1,\beta_2,\beta_3)=(\alpha_1*\beta_1,\beta_2,\beta_3),\\
(\alpha_1,\alpha_2,\alpha_3)\triangleright(\beta_1,\beta_2,\beta_3)=\bigl(\alpha_1\prec\beta_2,\alpha_2\star \beta_2^{\star-1},\alpha_3\bigr).
\end{gather*}
So $\Omega=(\eas(\Omega_1,*)\rtimes_\prec \eas'(\Omega_2,\star))\times \eas(\Omega_3)$.

In the particular case of $\eas(\Omega_1,*) \rtimes_\prec \eas'(\Omega_2,\star)$, \eqref{eq27} is trivial, and \eqref{eq28}, \eqref{eq29} can be reformulated in this way: for any $\alpha_1,\beta_1\in \Omega_1$, $\beta_2,\gamma_2\in \Omega_2$,
\begin{gather*}
(\alpha_1*\beta_1)\prec\gamma_2 =(\alpha_1\prec \gamma_2)*(\beta_1\prec \gamma_2),\qquad
(\alpha_1\prec \gamma_2)\prec \beta_2 =\alpha_1\prec (\beta_2\star \gamma_2).
\end{gather*}
As $\psi_{e_{\Omega_2}}=\operatorname{Id}_{\Omega_1}$, the following map is a left action of $(\Omega_2,\star)$ on $(\Omega_1,*)$ by group automorphisms:
\[\succ\colon \ \begin{cases}
\Omega_2\times \Omega_1 \longrightarrow \Omega_1,\\
(\beta_2,\alpha_1) \longrightarrow \beta_2\succ \alpha_1=\alpha_1\prec \beta_2.
\end{cases}\]
The formulas for the products in $\Omega$ are then immediate.
\end{proof}

\begin{Remark}
Consequently, we have a semi-direct product of groups $(\Omega_1,*)\rtimes_\succ (\Omega_2,\star)$.
\end{Remark}

Inverting the corresponding maps $\phi$, we obtain the following corollary.

\begin{Corollary}\label{cor3.15}
Let $\Omega$ be a finite non-degenerate dual CEDS. There exists an abelian group $(\Omega_1,*)$, a group $(\Omega_2,\star)$, a right action $\prec\colon\Omega_1\times \Omega_2\longrightarrow \Omega_1$ of $(\Omega_2,\star)$ on $(\Omega_1,*)$ by group automorphisms,
 and a nonempty set $\Omega_3$ such that $\Omega$ is of the form
\[(\eas(\Omega_2,\star) \ltimes_\prec \bigl(\eas'(\Omega_1,*)\bigr)\times \eas(\Omega_3),\]
with the products given by
\begin{gather*}
(\alpha_2,\alpha_1,\alpha_3)\rightarrow(\beta_2,\beta_1,\beta_3) =(\alpha_2\star\beta_2,\beta_1\prec\alpha_2,\beta_3),\\
(\alpha_2,\alpha_1,\alpha_3)\triangleright(\beta_2,\beta_1,\beta_3) =\bigl(\alpha_2,\alpha_1*\bigl(\beta_1^{-1}\prec \alpha_2^{-1}\bigr),\alpha_3\bigr).
\end{gather*} \end{Corollary}

\begin{Remark}
The inverse dual CEDS of the CEDS $(\eas(\Omega_1,*) \rtimes_\succ \eas'(\Omega_2,\star))\times \eas(\Omega_3)$ is $(\eas(\Omega_2,\star^{\rm op}) \ltimes_{\succ^{\rm op}} (\eas'(\Omega_1,*))\times \eas(\Omega_3)$.
\end{Remark}

\section{Linear extended associative semigroups}

\subsection{Definitions and examples}

The notions of $\ell$EAS, $\ell$CEDS and dual $\ell$CEDS are introduced in \cite[Definition 1.5]{Foissy48}, as a linear version of Lemma \ref{lem2.6}.

\begin{Definition}\label{defi4.1}
Let $A$ be a vector space and let $\Phi\colon A\otimes A\longrightarrow A\otimes A$ be a linear map.
\begin{enumerate}\itemsep=0pt
\item[(1)] We shall say that $(A,\Phi)$ is a linear extended associative semigroup (briefly, $\ell$EAS) if
\begin{align}
\label{eq30}
(\operatorname{Id} \otimes \Phi)\circ (\Phi\otimes \operatorname{Id})\circ (\operatorname{Id} \otimes \Phi)&=(\Phi \otimes \operatorname{Id})\circ (\operatorname{Id} \otimes \tau)\circ (\Phi\otimes \operatorname{Id}).
\end{align}
\item[(2)] We shall say that $(A,\Phi)$ is a linear commutative extended diassociative semigroup (briefly, $\ell$CEDS) if
\begin{gather}
\tag{\ref{eq30}}(\operatorname{Id} \otimes \Phi)\circ (\Phi\otimes \operatorname{Id})\circ (\operatorname{Id} \otimes \Phi) =(\Phi \otimes \operatorname{Id})\circ (\operatorname{Id} \otimes \tau)\circ (\Phi\otimes \operatorname{Id}),\\
(\operatorname{Id} \otimes \Phi)\circ (\operatorname{Id} \otimes \tau) \circ (\tau \otimes \operatorname{Id})\circ (\Phi \otimes \operatorname{Id})\nonumber\\
\qquad =(\tau \otimes \operatorname{Id})\circ (\Phi \otimes \operatorname{Id})\circ (\operatorname{Id} \otimes \Phi)\circ (\operatorname{Id} \otimes \tau).\label{eq31}
\end{gather}
\item[(3)] We shall say that $(A,\Phi)$ is a linear dual commutative extended diassociative semigroup (briefly, dual $\ell$CEDS) if
\begin{gather}
\tag{\ref{eq30}}(\operatorname{Id} \otimes \Phi)\circ (\Phi\otimes \operatorname{Id})\circ (\operatorname{Id} \otimes \Phi)=(\Phi \otimes \operatorname{Id})\circ (\operatorname{Id} \otimes \tau)\circ (\Phi\otimes \operatorname{Id}),\\
(\Phi \otimes \operatorname{Id})\circ (\tau \otimes \operatorname{Id})\circ (\operatorname{Id} \otimes \tau)\circ (\operatorname{Id} \otimes \Phi)\nonumber\\
\qquad=(\operatorname{Id} \otimes \tau)\circ (\operatorname{Id} \otimes \Phi)\circ (\Phi \otimes \operatorname{Id})\circ (\tau \otimes \operatorname{Id}).\label{eq32}
\end{gather}\end{enumerate}
If $(A,\Phi)$ is an $\ell$EAS (resp.\ an $\ell$CEDS or a dual $\ell$CEDS), we shall say that it is non-degenerate if $\Phi$ is bijective.
\end{Definition}

Note that, by definition, $\ell$CEDS and dual $\ell$CEDS are $\ell$EAS.

\begin{Example}\label{ex4.1} \quad
\begin{enumerate}\itemsep=0pt
\item[(1)] Let $(\Omega,\rightarrow,\triangleright)$ be an EAS (resp.\ a CEDS, a dual CEDS). Let $A=\mathbb{K} \Omega$ be the vector space generated by $\Omega$. We define
\[\Phi\colon \ \begin{cases}
A\otimes A \longrightarrow A\otimes A,\\
a\otimes b \longrightarrow (a\rightarrow b)\otimes (a\triangleright b).
\end{cases}\]
Then $(A,\Phi)$ is an $\ell$EAS (resp.\ an $\ell$CEDS, a dual $\ell$CEDS), which we call the linearization of $(\Omega,\rightarrow,\triangleright)$. It is a non-degenerate $\ell$EAS if and only if $\Omega$ is a non-degenerate EAS.
\item[(2)] Not all the $\ell$EAS are of the form $\mathbb{K}\Omega$. For example, if $A$ is a two-dimensional space with basis $(x,y)$, the maps given by the following matrices in the basis $(x\otimes x,x\otimes y,y\otimes x,y\otimes y)$ are $\ell$EAS
\begin{alignat*}{4}
&M_1 =\begin{pmatrix}0&0&1&0\\0&0&0&0\\0&0&0&0\\0&0&0&0\end{pmatrix},\qquad&&
M_2 =\begin{pmatrix}0&0&0&0\\0&0&a&0\\0&0&0&0\\0&0&0&0\end{pmatrix},\qquad&&
M_3=\begin{pmatrix}1&0&0&0\\0&0&0&0\\0&0&0&0\\0&0&0&0\end{pmatrix},&\\
&M_4=\begin{pmatrix}1&0&0&0\\0&0&1&0\\0&0&0&0\\0&0&0&0\end{pmatrix},\qquad&&
M_5=\begin{pmatrix}1&0&0&0\\0&0&0&0\\0&1&0&0\\0&0&0&0\end{pmatrix},\qquad&&
M_6=\begin{pmatrix}1&0&0&0\\0&0&0&0\\0&0&1&0\\0&0&0&0\end{pmatrix},&\\
&M_7=\begin{pmatrix}1&0&0&0\\0&0&0&0\\0&0&0&0\\0&0&1&0\end{pmatrix},\qquad&&
M_8=\begin{pmatrix}1&0&0&0\\0&0&1&0\\0&0&1&0\\0&0&0&0\end{pmatrix},\qquad&&
M_9=\begin{pmatrix}1&0&0&0\\0&0&0&0\\0&1&0&0\\0&0&1&0\end{pmatrix},&\\
&M_{10}=\begin{pmatrix}1&0&0&0\\0&0&0&0\\0&1&1&0\\0&0&0&0\end{pmatrix},\qquad&&
M_{11}=\begin{pmatrix}1&0&0&0\\0&0&1&0\\0&0&1&0\\0&0&1&0\end{pmatrix},\qquad&&
M_{12}=\begin{pmatrix}1&0&0&0\\0&0&0&0\\0&1&0&0\\0&1&-1&0\end{pmatrix},&\\
&M_{13}=\begin{pmatrix}1&0&0&0\\0&0&0&0\\0&0&0&0\\0&0&0&1\end{pmatrix},\qquad&&
M_{14}=\begin{pmatrix}1&0&0&0\\0&0&0&0\\0&1&0&0\\0&0&0&1\end{pmatrix},\qquad&&
M_{15}=\begin{pmatrix}1&1&0&0\\0&0&0&0\\0&0&0&0\\0&0&1&1\end{pmatrix},&\\
&M_{16}=\begin{pmatrix}1&0&0&0\\0&0&0&0\\0&1&0&0\\0&0&1&1\end{pmatrix},\qquad&&
M_{17}=\begin{pmatrix}1&0&1&0\\0&0&-1&0\\0&1&-1&0\\0&0&2&1\end{pmatrix},\qquad&&
M_{18}=\begin{pmatrix}1&0&0&0\\0&0&1&0\\0&1&0&0\\0&0&0&1\end{pmatrix},&
\end{alignat*}
where $a$ is a scalar. Moreover,
\begin{itemize}\itemsep=0pt
\item The $\ell$CEDS in this list are the $M_i$'s with
\[i\in \{1,2,3,4,5,9,10,13,14,16,17,18\}.\]
\item The dual $\ell$CEDS in this list are the $M_i$'s with
\[i\in \{1,2,3,4,5,7,8,9,11,13,14,15,\allowbreak 16,17,18\}.\]
\end{itemize}
These $\ell$EAS are in fact the EAS of dimension 2 which have a basis of special vectors, see Definition \ref{defi4.3}.
\end{enumerate}\end{Example}

\begin{Notation}
Let $(A,\Phi)$ be an $\ell$EAS. We use Sweedler's-like notation
\[
\Phi(a\otimes b)=\sum a'\rightarrow b' \otimes a''\triangleright b''.
\]
Note that the operations $\rightarrow$ and $\triangleright$ do not necessarily exist, nor the coproducts $a'\otimes a''$ or $b'\otimes b''$.
With this notation, \eqref{eq30} can be rewritten as
\begin{gather}
\sum\sum\sum a' \rightarrow \bigl(b' \rightarrow c'\bigr)'\otimes \bigl(a'' \triangleright \bigl(b' \rightarrow c'\bigr)''\bigr)'\rightarrow \bigl(b'' \triangleright c''\bigr)' \otimes \bigl(a'' \triangleright \bigl(b' \rightarrow c'\bigr)\bigr)''\triangleright \bigl(b'' \triangleright c''\bigr)''\nonumber\\
 \qquad=\sum\sum \bigl(a' \rightarrow b'\bigr)'\rightarrow c' \otimes \bigl(a' \rightarrow b'\bigr)''\triangleright c'' \otimes a'' \triangleright b''.\tag{\ref{eq30}$'$}
\end{gather}
Similarly, \eqref{eq31} and \eqref{eq32} are rewritten as
\begin{gather}
 \sum\sum a'' \triangleright \bigl(c' \rightarrow b'\bigr)'' \otimes a' \rightarrow \bigl(c'' \rightarrow b''\bigr)'\otimes c'' \triangleright b''\nonumber\\
\qquad=\sum\sum a'' \triangleright b'' \otimes c' \rightarrow \bigl(a' \rightarrow b'\bigr)'\otimes c'' \triangleright \bigl(a' \rightarrow b'\bigr)'',\tag{\ref{eq31}$'$}\\
\nonumber\\
 \sum\sum \bigl(b'' \triangleright c''\bigr)'\rightarrow a'\otimes \bigl(b'' \triangleright c''\bigr)''\triangleright a''\otimes b' \rightarrow c'\nonumber\\
\qquad=\sum\sum b' \rightarrow a' \otimes \bigl(b'' \triangleright a''\bigr)''\triangleright c''\otimes \bigl(b'' \triangleright a''\bigr)'\rightarrow c'.\tag{\ref{eq32}$'$}
\end{gather}\end{Notation}

By transposition of \eqref{eq30}, \eqref{eq31} and \eqref{eq32}, we have the following.

\begin{Proposition}\label{prop4.2}
Let $V$ be a finite-dimensional space and
\[
\Phi\colon \ V\otimes V\longrightarrow V\otimes V
\]
 be a linear map. We consider
\[
\Phi^*\colon \ V^*\otimes V^*=(V\otimes V)^*\longrightarrow(V\otimes V)^*=V^*\otimes V^*.
\]
Then $(V,\Phi)$ is an $\ell$EAS $($resp.\ an $\ell$CEDS, a dual $\ell$CEDS$)$ if and only if $(V^*,\Phi^*)$ is an $\ell$EAS $($resp.\ a~dual $\ell$CEDS, an $\ell$CEDS$)$.
\end{Proposition}

\begin{Example}\quad
\begin{enumerate}\itemsep=0pt
\item[(1)] As their matrices are symmetric, the $\ell$EAS $M_3$, $M_6$, $M_{13}$ and $M_{18}$ are self-dual, through the pairing which matrix in the basis $(x,y)$ is
$\bigl(\begin{smallmatrix}
1&0\\0&1
\end{smallmatrix}\bigr)$.
With the same pairing, the dual of $M_4$ is $M_5$ and the dual of $M_8$ is $M_{10}$. The $\ell$EAS $M_2$ and $M_{14}$ are also self-dual, through the pairing which matrix in the basis~$(x,y)$ is
$\bigl(\begin{smallmatrix}
0&1\\1&0
\end{smallmatrix}\bigr)$.
The $\ell$EAS $M_{16}$ and $M_{17}$ are self-dual,\footnote{For $M_{17}$, this holds if the characteristic of the base field is not 2.} through the pairings which matrix in the basis $(x,y)$ are respectively
$\bigl(\begin{smallmatrix}
1&1\\1&0
\end{smallmatrix}\bigr)$, $\bigl(\begin{smallmatrix}
2&2\\2&1
\end{smallmatrix}\bigr)$.
The duals of $M_1$, $M_7$, $M_9$, $M_{11}$, $M_{12}$ and $M_{15}$ are not isomorphic to any $M_i$'s.
\item[(2)] Let $\Omega$ be a finite EAS and $A=\mathbb{K}\Omega$ be the associated $\ell$EAS. The dual $A^*$ is identified with the space $\mathbb{K}^\Omega$ of maps from $\Omega$ to $\mathbb{K}$, with the dual basis $(\delta_\alpha)_{\alpha \in \Omega}$ of the basis~$\Omega$ of~$\mathbb{K}\Omega$.
Then, for any $\alpha,\beta \in \Omega$,
\[
\Phi^*(\delta_\alpha\otimes \delta_\beta)=\sum_{(\gamma,\delta)\in \phi^{-1}(\alpha,\beta)} \delta_\gamma \otimes \delta_\delta.
\]
This is usually not the linearization of an EAS, except if $\Omega$ is non-degenerate: in this case, we recover the linearization of $(\Omega,\curvearrowright,\blacktriangleright)$ of Proposition \ref{prop2.7}.
\end{enumerate}\end{Example}

\subsection{Special vectors, left units and counits}

\begin{Definition}\label{defi4.3}
Let $(A,\Phi)$ be an $\ell$EAS.
\begin{enumerate}\itemsep=0pt
\item[(1)] Let $a\in A$. We shall say that $a$ is a left unit of $(A,\Phi)$ if for any $b\in A$,
\[
\Phi(a\otimes b)=b\otimes a.
\]
\item[(2)] Let $f \in A^*$. We shall say that $f$ is a left counit of $(A,\Phi)$ if
\[
(f \otimes \operatorname{Id})\circ \Phi=\operatorname{Id} \otimes f.
\]
\item[(3)] Let $a\in A$ and $\lambda \in \mathbb{K}$. We shall say that $a$ is a special vector of $(A,\Phi)$ of eigenvalue $\lambda$ if
 \[
 \Phi(a\otimes a)=\lambda a\otimes a.
 \]
\end{enumerate}
\end{Definition}

\begin{Remark}
Let $(A,\Phi)$ be an $\ell$EAS.
\begin{enumerate}\itemsep=0pt
\item[(1)] Any left unit of $(A,\Phi)$ is a special vector of eigenvalue $1$.
\item[(2)] If $A$ is finite-dimensional, its left counits are the left units of $(A^*,\Phi^*)$.
\item[(3)] The set of left units is a subspace of $A$ and the set of left counits a subspace of $A^*$. The set of special vectors of a given eigenvalue is generally not a subspace of $A$.
\end{enumerate}
\end{Remark}

\begin{Lemma}
Let $(A,\Phi)$ be an $\ell$EAS and $a\in A$ be a nonzero special vector of $(A,\Phi)$. Then its eigenvalue $\lambda$ is $0$ or $1$.
\end{Lemma}

\begin{proof}
Let us apply \eqref{eq30} to $a\otimes a\otimes a$. This gives
\[
\lambda^3 a\otimes a\otimes a=\lambda^2 a\otimes a\otimes a.
\] As $a\neq 0$, $\lambda=0$ or $1$.\end{proof}

\begin{Example}\quad
\begin{enumerate}\itemsep=0pt
\item[(1)] Let us give special vectors, left units and left counits for the thirteen $\ell$EAS associated to the thirteen EAS of cardinality 2. In each case, we give a basis of the spaces of left units and left counits; $\lambda$, $\mu$ and $\nu$ are scalars.
The dual basis of the basis $(X,Y)$ of $\mathbb{K}\Omega$ is denoted by $(X^*,Y^*)$.
$$\renewcommand{\arraystretch}{1.2} \begin{array}{|c||c|c|c|c|}
\hline \text{Case}&\text{Special vectors}&\text{Special vectors}&\text{Left units}&\text{Left conits}\\
&\text{of eigenvalue 1}&\text{of eigenvalue 0}&&\\
\hline \hline \mathbf{A1}&\lambda X&\nu(X-Y)&\varnothing&\varnothing\\
\hline \mathbf{A2}&\lambda X&\nu(X-Y)&\varnothing&(X^*+Y^*)\\
\hline \mathbf{C1}&\lambda X&0&\varnothing&\varnothing\\
\hline \mathbf{C3}&\lambda X,\, \mu Y&0&(Y)&(X^*+Y^*)\\
\hline \mathbf{C5}&\mu Y&0&(Y)&\varnothing\\
\hline \mathbf{C6}&\lambda X&0&\varnothing&\varnothing\\
\hline \mathbf{E'1}\text{--}\mathbf{E'2}&\lambda X&\nu(X-Y)&\varnothing&\varnothing\\
\hline \mathbf{E'3}&\lambda X,\, \mu Y&\nu(X-Y)&\varnothing&(X^*+Y^*)\\
\hline \mathbf{F1}&\lambda X&\nu(X-Y)&(X)&\varnothing\\
\hline \mathbf{F3}&\lambda X+\mu Y&0&(X,Y)&(X^*,Y^*)\\
\hline \mathbf{F4}&\lambda X,\, \nu(X+Y)&0&(X+Y)&(X^*)\\
\hline \mathbf{H1}&\lambda X&0&(X)&\varnothing\\
\hline \mathbf{H2}&\lambda X,\, \nu(X+Y)&0&(X)&(X^*+Y^*)\\
\hline \end{array}$$
Some of them have a basis of special vectors: let us determine their matrices in such a~basis. We recover in this way some matrices of Example~\ref{ex4.1}:
\begin{itemize}\itemsep=0pt
\item For $\mathbf{A1}$, in the basis $(X,X-Y)$, we obtain $M_3$.
\item For $\mathbf{A2}$, in the basis $(X,X-Y)$, and for $\mathbf{F1}$, in the basis $(X-Y,X)$, we obtain $M_4$.
\item For $\mathbf{C3}$, in the basis $(Y,X)$, we obtain $M_{16}$.
\item For $\mathbf{E'1}\text{--}\mathbf{E'2}$, in the basis $(X,X-Y)$, we obtain $M_6$.
\item For $\mathbf{E'3}$, in the basis $(X,Y-X)$, we obtain $M_{11}$.
\item For $\mathbf{F3}$, in the basis $(X,Y)$, we obtain $M_{18}$.
\item For $\mathbf{F4}$ and $\mathbf{H2}$, in the basis $(X-Y,X)$, we obtain $M_{17}$.
\end{itemize}
Hence, the $\ell$EAS associated to $\mathbf{A2}$ and $\mathbf{F1}$ are isomorphic, whereas the EAS $\mathbf{A2}$ and $\mathbf{F1}$ are not. As similar situation holds for $\mathbf{F4}$ and $\mathbf{H2}$.

\item[(2)] It is possible to show that any 2-dimensional $\ell$EAS with a basis of special vectors is isomorphic to one of the eighteen cases of Example \ref{ex4.1}.
For all of them, let us give special vectors, left units and left counits for the eighteen cases of Example~\ref{ex4.1}.
In each case, we give a basis of the spaces of left units and left counits; $\lambda$, $\mu$ and $\nu$ are scalars. For~$M_2$, we~assume that the parameter $a$ is nonzero (otherwise, $\Phi=0$).
$$\renewcommand{\arraystretch}{1.2} \begin{array}{|c||c|c|c|c|}
\hline \text{Case}&\text{Special vectors}&\text{Special vectors}&\text{Left units}&\text{Left counits}\\
&\text{of eigenvalue 1}&\text{of eigenvalue 0}&&\\
\hline\hline M_1&0&\lambda x, \mu y&\varnothing&\varnothing\\
\hline M_2&0&\lambda x, \mu y&\varnothing&\varnothing\\
\hline M_3&\lambda x&\mu y&\varnothing&\varnothing\\
\hline M_4&\lambda x&\mu y&\varnothing& (x^*)\\
\hline M_5&\lambda x&\mu y&(x)&\varnothing\\
\hline M_6&\lambda x&\mu y&\varnothing&\varnothing\\
\hline M_7&\lambda x&\mu y&\varnothing&\varnothing\\
\hline M_8&\lambda x&\mu y&\varnothing& (x^*)\\
\hline M_9&\lambda x&\mu y&(x)&\varnothing\\
\hline M_{10}&\lambda x&\mu y&(x)&\varnothing\\
\hline M_{11}&\lambda x,\nu(x+y)&\mu y&\varnothing& (x^*)\\
\hline M_{12}&\lambda x&\mu y&\varnothing&\varnothing\\
\hline M_{13}&\lambda x, \mu y&0&\varnothing&\varnothing\\
\hline M_{14}&\lambda x, \mu y&0&(x)&(y^*)\\
\hline M_{15}&\lambda x, \mu y&\nu(x-y)&\varnothing&(x^*+y^*)\\
\hline M_{16}&\lambda x, \mu y&0&(x)&(x^*+y^*)\\
\hline M_{17}&\lambda x, \mu y&0&(x)&(x^*+y^*)\\
\hline M_{18}&\lambda x+\mu y&0&(x,y)&(x^*,y^*)\\
\hline \end{array}$$
Among them, $M_{11}$ has three lines of special vectors. In the basis $(x+y,x)$ its matrix is~$M_{15}$, so $M_{11}$ and $M_{15}$ are isomorphic.
\end{enumerate}\end{Example}

\subsection{Left units and counits of finite non-degenerate CEDS}

\begin{Proposition}\label{prop4.5}
Let $(\Omega_1,*)$ be an abelian finite group, $(\Omega_2,\star)$ be a finite group, and $\Omega_3$ be a~finite set. We denote by $e_1$ and $e_2$ the units of $\Omega_1$ and $\Omega_2$.
\begin{enumerate}\itemsep=0pt
\item[$(1)$] Let $(A,\Phi)$ be the linearization of the CEDS
$(\eas(\Omega_1,*) \rtimes_\succ \eas'(\Omega_2,\star))\times \eas(\Omega_3)$
of Theorem {\rm\ref{theo3.14}}.
\begin{enumerate}\itemsep=0pt
\item[$(a)$] The special vectors of eigenvalue $1$ of $(A,\Phi)$ are the vectors of the form
\[a=\sum_{(\alpha_1,\alpha_2,\alpha_3)\in H_1\times H_2\times \Omega_3}g(\alpha_3)(\alpha_1,\alpha_2,\alpha_3),\]
where $H_1$ is a subgroup of $\Omega_1$, $H_2$ is a subgroup of $\Omega_2$, such that $H_2\succ H_1\subseteq H_1$, and~${g\colon \Omega_3\longrightarrow \mathbb{K}}$ is a map.
\item[$(b)$] The left units of $(A,\Phi)$ are the vectors of the form
\[a=\sum_{(\alpha_2,\alpha_3)\in \Omega_2\times \Omega_3}g(\alpha_3)(e_1,\alpha_2,\alpha_3),\]
where $g\colon \Omega_3\longrightarrow \mathbb{K}$ is a map.
\item[$(c)$] The left counits of $(A,\Phi)$ are the linear forms $f$ such that for any $(\alpha_1,\alpha_2,\alpha_3)\in \Omega$,
$f(\alpha_1,\alpha_2,\alpha_3)=\delta_{\alpha_2,e_2}g(\alpha_3)$,
where $g\colon \Omega_3\longrightarrow \mathbb{K}$ is a map.
\end{enumerate}
\item[$(2)$] Let $(A,\Phi)$ be the linearization of the CEDS
$(\eas(\Omega_2,\star) \ltimes_\prec \eas'(\Omega_1,*))\times \eas(\Omega_3)$
of Corollary {\rm\ref{cor3.15}}.
\begin{enumerate}\itemsep=0pt
\item[$(a)$] The special vectors of eigenvalue 1 of $(A,\Phi)$ are the vectors of the form
\[a=\sum_{(\alpha_1,\alpha_2,\alpha_3)\in H_1\times H_2\times \Omega_3}g(\alpha_3)(\alpha_2,\alpha_1,\alpha_3),\]
where $H_1$ is a subgroup of $\Omega_1$, $H_2$ is a subgroup of $\Omega_2$, such that $H_1\prec H_2\subseteq H_1$, and~${g\colon \Omega_3\longrightarrow \mathbb{K}}$ is a map.
\item[$(b)$] The left units of $(A,\Phi)$ are the vectors of the form
\[a=\sum_{(\alpha_1,\alpha_3)\in \Omega_1\times \Omega_3}g(\alpha_3)(e_2,\alpha_1,\alpha_3),\]
where $g\colon \Omega_3\longrightarrow \mathbb{K}$ is a map.
\item[$(c)$] The left counits of $(A,\Phi)$ are the linear forms $f$ such that for any $(\alpha_1,\alpha_2,\alpha_3)\in \Omega$,
$f(\alpha_1,\alpha_2,\alpha_3)=\delta_{\alpha_1,e_1}g(\alpha_3)$,
where $g\colon \Omega_3\longrightarrow \mathbb{K}$ is a map.
\end{enumerate} \end{enumerate}\end{Proposition}

\begin{proof}
(1)(a) Let $a$ be a nonzero vector of $A$, which we write as
\[a=\sum_{(\alpha_1,\alpha_2,\alpha_3)\in \Omega}a_{(\alpha_1,\alpha_2,\alpha_3)}(\alpha_1,\alpha_2,\alpha_3).\]
Then $a$ is a special vector of eigenvalue~1 if and only if for any $(\alpha_1,\alpha_2,\alpha_3)$, $(\beta_1,\beta_2,\beta_3)\in \Omega$,
\begin{align}
\label{eq33}a_{(\alpha_1,\alpha_2,\alpha_3)}a_{(\beta_1,\beta_2,\beta_3)}&=a_{(\alpha_2^{-1}\succ \beta_1,\beta_2\star \alpha_2,\beta_3)}a_{((\alpha_2^{-1}\succ \beta_1^{-1})*\alpha_1,\alpha_2,\alpha_3)}.
\end{align}
We put
\begin{gather*}
 \forall (\alpha_2,\alpha_3)\in \Omega_2\times \Omega_3,\qquad H_1(\alpha_2,\alpha_3)=\{\alpha_1\in \Omega_1, a_{(\alpha_1,\alpha_2,\alpha_3)}\neq 0\},\\
 \forall \alpha_3\in \Omega_3,\qquad H_2(\alpha_3)=\{\alpha_2\in \Omega_2, H_1(\alpha_2,\alpha_3)\neq \varnothing\},\qquad
 H_3=\{\alpha_3\in \Omega_3, H_2(\alpha_3)\neq \varnothing\}.
\end{gather*}
We shall also consider the map
\[g\colon \ \begin{cases}
\Omega_3 \longrightarrow \mathbb{K},\\
\alpha_3 \longrightarrow g(\alpha_3)=a_{e_1,e_2,\alpha_3}.
\end{cases}\]

Let us first prove that if $\alpha_3\in H_3$, then $H_2(\alpha_3)$ is a subgroup of $\Omega_2$. Let $\alpha_2,\beta_2\in H_2(\alpha_3)$ (which is nonempty as $\alpha_3\in H_3$).
Let $\alpha_1,\beta_1\in \Omega_1$, such that $a_{(\alpha_1,\alpha_2,\alpha_3)}\neq 0$ and $a_{(\beta_1,\beta_2,\alpha_3)}\neq 0$. By \eqref{eq33},
\smash{$a_{(\alpha_2^{-1}\succ \beta_1,\beta_2\star \alpha_2,\alpha_3)}\neq 0$},
so $\beta_2\star \alpha_2\in H_2(\alpha_3)$. As $\Omega_2$ is finite, $H_2(\alpha_3)$ is a subgroup of~$\Omega_2$.

Let us prove that if $\alpha_3\in H_3$, then $H_1(e_2,\alpha_3)$ is a subgroup of $\Omega_1$ and, moreover, for any $\alpha_1\in H_1(e_2,\alpha_3)$, $a_{(\alpha_1,e_2,\alpha_3)}=g(\alpha_3)$.
As $H_2(\alpha_3)$ is a subgroup of $\Omega_2$, it contains $e_2$, so $H_1(e_2,\alpha_3)\neq \varnothing$. Let $\alpha_1,\beta_1\in H_1(e_2,\alpha_3)$. By \eqref{eq33},
\[a_{(\alpha_1,e_2,\alpha_3)}a_{(\beta_1,e_2,\alpha_3)}=a_{(\beta_1,e_2,\alpha_3)}a_{(\beta_1^{-1}*\alpha_1,e_2,\alpha_3)}\neq 0.\]
Hence, $\beta_1^{-1}*\alpha_1\in H_1(e_2,\alpha_3)$. Taking $\alpha_1=\beta_1$, we obtain{\samepage
\[a_{(\alpha_1,e_2,\alpha_3)}a_{(\alpha_1,e_2,\alpha_3)}=a_{(\alpha_1,e_2,\alpha_3)}a_{(e_1,e_2,\alpha_3)}\neq 0,\]
so $a_{(\alpha_1,e_2,\alpha_3)}=a_{(e_1,e_2,\alpha_3)}=g(\alpha_3)$.}

Let us prove that if $\alpha_3,\beta_3\in H_3$ and $\beta_2\in H_2(\beta_3)$, then $H_1(\beta_2,\beta_3)\subseteq H_1(e_2,\alpha_3)$. Let $\beta_1\in H_1(\beta_2,\beta_3)$.
Then $a_{(\beta_1,\beta_2,\beta_3)}\neq 0$. As $H_1(e_1,\beta_2)$ is a subgroup of $\Omega_1$, it contains $e_1$, so~${a_{(e_1,e_2,\alpha_3)}\neq 0}$. By \eqref{eq33},
\[a_{(e_1,e_2,\alpha_3)}a_{(\beta_1,\beta_2,\beta_3)}=g(\alpha_3)a_{(\beta_1,\beta_2,\beta_3)}=a_{(\beta_1,\beta_2,\beta_3)}a_{(\beta_1^{-1},e_2,\alpha_3)}\neq 0,\]
so $\beta_1^{-1}\in H_1(e_2,\alpha_3)$. As this is a subgroup of $\Omega_1$, $\beta_1\in H_1(e_2,\alpha_3)$.
As a consequence, for $\beta_2=e_2$, we obtain by symmetry that for any $\alpha_3,\beta_3\in H_3$, $H_1(e_2,\alpha_3)=H_1(e_2,\beta_3)$.
Therefore, there exists a subgroup $H_1$ of $\Omega_1$ such that for any $\alpha_3\in \Omega_3$, $H_1(e_2,\alpha_3)=H_1$.

Let us prove that for any $\alpha_3\in H_3$, $H_2(\alpha_3)\succ H_1\subseteq H_1$. Let $\beta_1\in H_1=H_1(e_2,\alpha_3)$, then $a_{(\beta_1,e_2,\alpha_3)}\neq 0$. Let $\alpha'_2\in H_2(\alpha_3)$.
We put $\alpha_2=\alpha_2'^{-1}\in H_2(\alpha_3)$. There exists $\alpha_1\in H_1(\alpha_2,\alpha_3)$, such that $a_{(\alpha_1,\alpha_2,\alpha_3)}\neq 0$. By \eqref{eq33},
 \[a_{(\alpha_1,\alpha_2,\alpha_3)}a_{(\beta_1,e_2,\alpha_3)}=a_{(\beta_1,\alpha_2,\alpha_3)}a_{((\alpha_2^{-1}\succ \beta_1^{-1})*\alpha_1,e_2,\alpha_3)}\neq 0,\]
so $\bigl(\alpha_2^{-1}\succ \beta_1^{-1}\bigr)*\alpha_1\in H_1(e_2,\alpha_3)=H_1$. Moreover, as $H_1(\alpha_2,\alpha_2)\subseteq H_1$,
\[\alpha_2^{-1}\succ \beta_1^{-1}=\alpha_2'\succ \beta_1^{-1}\in H_1.\]
Its inverse $ \alpha'_2\succ \beta_1$ is also an element of $H_1$, so $H_2(\alpha_3)\succ H_1\subseteq H_1$.

Let us prove that for any $\alpha_3\in H_3$, for any $\alpha_2\in H_2(\alpha_3)$, $H_1(\alpha_2,\alpha_3)=H_1$. Let $\alpha_1\in H_1=H_1(e_2,\alpha_3)$ and $\beta_1\in H_1(\alpha_2,\alpha_3)$. By \eqref{eq33},
\[a_{(\alpha_1,e_2,\alpha_3)}a_{(\beta_1,\alpha_2,\alpha_3)}=a_{(\beta_1,\alpha_2,\alpha_3)}a_{(\beta_1^{-1}*\alpha_1,\alpha_2,\alpha_3)}\neq 0,\]
so $\beta_1^{-1}*\alpha_1\in H_1(\alpha_2,\alpha_3)$. We obtain an injective map
\[\begin{cases}
H_1 \longrightarrow H_1(\alpha_2,\alpha_3),\\
\alpha_1 \longrightarrow \beta_1^{-1}*\alpha_1.
\end{cases}\]
Hence, $|H_1|\leqslant |H_1(\alpha_2,\alpha_3)|$. We already proved that $H_1(\alpha_2,\alpha_3)\subseteq H_1$, so $H_1=H_1(\alpha_2,\alpha_3)$.

We now prove that there exists a subgroup $H_2$ of $\Omega_2$ such that for any $\alpha_3\in H_3$, $H_2(\alpha_3)=H_2$, and that, moreover, for any $\alpha_2\in H_2$, for any $\alpha_3\in H_3$, $a_{(e_1,\alpha_2,\alpha_3)}=g(\alpha_3)$.
Let $\alpha_3,\beta_3\in H_3$. Let $\alpha_2\in H_2(\alpha_3)$. As $e_1\in H_1(\alpha_2,\alpha_3)=H_1$, $a_{(e_1,\alpha_2,\alpha_3)}\neq 0$. As $e_2\in H_2(\alpha_3)$, $a_{(e_1,e_2,\beta_3)}\neq 0$. By~\eqref{eq33},
\[a_{(e_1,\alpha_2,\alpha_3)}a_{(e_1,e_2,\beta_3)}=a_{(e_1,\alpha_2,\beta_3)}a_{(e_1,\alpha_2\alpha_3)}\neq 0,\]
so $\alpha_2\in H_2(\beta_3)$. We proved that $H_2(\alpha_3)\subseteq H_2(\beta_3)$: by symmetry, $H_2(\alpha_3)=H_2(\beta_3)$, which prove the existence of $H_2$.
Moreover, as $a_{(e_1,\alpha_2,\alpha_3)}\neq 0$ and $a_{(e_1,e_2,\beta_3)}=g(\beta_3)$, we obtain that $a_{(e_1,\alpha_2,\beta_3)}=g(\beta_3)$.

We proved that for any $(\alpha_1,\alpha_2,\alpha_3)\in \Omega$,
\[a_{(\alpha_1,\alpha_2,\alpha_3)}\neq 0\Longleftrightarrow (\alpha_1,\alpha_2,\alpha_3)\in H_1\times H_2\times H_3.\]
Let $\alpha_3,\beta_3\in H_3$, $\beta_1\in H_1$, $\alpha_2\in H_2$. We put $\alpha_1=\alpha_2^{-1}\succ \beta_1$. Then, by \eqref{eq33},
\[g(\alpha_3)g(\beta_3)=a_{(\beta_1,\alpha_2,\alpha_3)}a_{(e_1,\alpha_2,\beta_3)}=a_{(\beta_1,\alpha_2,\alpha_3)}g(\beta_3),\]
so $a_{(\beta_1,\alpha_2,\alpha_3)}=g(\alpha_3)$. We proved that $a$ has the announced form.

Conversely, if $a$ is of the announced form,
\begin{align*}
\Phi(a\otimes a)={}&\sum_{(\alpha_1,\alpha_2,\alpha_3),(\beta_1,\beta_2,\beta_3)\in H_1\times H_2\times \Omega_3}g(\alpha_3)g(\beta_3)(\alpha_1*\beta_1,\beta_2,\beta_3)\\
&\otimes \bigl(\beta_2\succ \alpha_2,\alpha_2*\beta_2^{-1},\alpha_3\bigr)\\
={}&\sum_{(\alpha_1,\alpha_2,\alpha_3),(\beta_1,\beta_2,\beta_3)\in H_1\times H_2\times \Omega_3}(\beta_1,\beta_2,\beta_3)\otimes (\alpha_1,\alpha_2,\alpha_3)
=a\otimes a.
\end{align*}

(2)(a) If $(A,\Phi)$ is the linearization of $(\eas(\Omega_2,\star) \ltimes_\prec \eas'(\Omega_1,*))\times \eas(\Omega_3)$, then $\bigl(A,\Phi^{-1}\bigr)$ is the linearization of $(\eas(\Omega_1,*) \rtimes_{\prec^{\rm op}} \eas'(\Omega_2,\star^{\rm op}))\times \eas(\Omega_3)$.
The result then comes from the observation that the special vectors of $(A,\Phi)$ and $\bigl(A,\Phi^{-1}\bigr)$ are the same.

(1)(b) Let $a$ be a left unit of $A$. Then it is a special vector, which we write as
\[a=\sum_{(\alpha_1,\alpha_1,\alpha_3)\in H_1\times H_2\times \Omega_3}g(\alpha_3)(\alpha_1,\alpha_1,\alpha_3).\]
For any $b=(\beta_1,\beta_2,\beta_3)\in \Omega$,
\begin{align*}
\Phi(a\otimes b)&=\sum_{(\alpha_1,\alpha_1,\alpha_3)\in H_1\times H_2\times \Omega_3}g(\alpha_3)(\alpha_1*\beta_1,\beta_2,\beta_3)\otimes \bigl(\beta_2\succ \alpha_1,\alpha_2*\beta_2^{-1},\alpha_3\bigr)\\
&=b\otimes a=\sum_{(\alpha_1,\alpha_1,\alpha_3)\in H_1\times H_2\times \Omega_3}g(\alpha_3)(\beta_1,\beta_2,\beta_3)\otimes (\alpha_1,\alpha_1,\alpha_3).
\end{align*}
Taking $\beta_1=e_1$, we obtain that for any $\alpha_1\in H_1$, $\alpha_1=e_1$, so $H_1=\{e_1\}$. Moreover, for any $\beta_2\in \Omega_2$,
\[\sum_{(\alpha_2,\alpha_3)\in H_2\times \Omega_3}g(\alpha_3)\bigl(\alpha_2\star \beta_2^{-1},\alpha_3\bigr)=\sum_{(\alpha_2,\alpha_3)\in H_2\times \Omega_3}g(\alpha_3)(\alpha_2,\alpha_3),\]
so for any $\alpha_2\in H_2$, $\alpha_2\star \beta_2^{-1}\in H_2$. In particular, for $\alpha_2=e_2$, $\beta_2^{-1}\in H_2$ and finally $\beta_2\in H_2$: $H_2=\Omega_2$. The converse application is immediate.

(2)(b) Similar proof.

(1)(c) and (2)(c) The left counits of $(A,\Phi)$ are the left units of $(A^*,\Phi^*)$, which is isomorphic to $\bigl(A,\Phi^{-1}\bigr)$.The result comes from(2)(b). and (1)(b). \end{proof}

\section[From bialgebras to lEAS]{From bialgebras to $\boldsymbol{\ell}$EAS}

We refer to \cite{Abe1980,Cartier2021,Sweedler1969} for classical results and notations on bialgebras and Hopf algebras.

\subsection[A functor from bialgebras to lEAS]{A functor from bialgebras to $\boldsymbol{\ell}$EAS}

\begin{Proposition} \label{prop5.1}
Let $(A,m,\Delta)$ be a bialgebra, not necessarily unitary nor counitary. For any $a,b\in A$, we define $\Phi\colon A\otimes A\longrightarrow A\otimes A$ by
\begin{align*}
\Phi(a\otimes b) =(m\otimes \operatorname{Id}_A)\circ (\operatorname{Id}_A\otimes \tau)\circ (\Delta\otimes \operatorname{Id}_A)(a\otimes b)=\sum a^{(1)}b\otimes a^{(2)},
\end{align*}
with Sweedler's notation $\Delta(a)=\sum a^{(1)}\otimes a^{(2)}$. Then $(A,\Phi)$ is an $\ell$EAS, denoted by $\leas(A,\allowbreak m,\Delta)$.
\end{Proposition}

\begin{proof}
For any $a,b,c\in A$,
\begin{align*}
(\operatorname{Id} \otimes \Phi)\circ(\Phi\otimes \operatorname{Id})\circ (\operatorname{Id} \otimes \Phi)(a\otimes b\otimes c)&=(\Phi\otimes \operatorname{Id})\circ (\operatorname{Id} \otimes \tau)\circ (\Phi \otimes \operatorname{Id})(a\otimes b\otimes c)\\
&=\sum \sum a^{(1)}b^{(1)}c\otimes a^{(2)}b^{(2)}\otimes a^{(3)}. \tag*{\qed}
\end{align*} \renewcommand{\qed}{}\end{proof}

\begin{Example}\quad\samepage
\begin{enumerate}\itemsep=0pt
\item[(1)] Let $(\Omega,\star)$ be a semigroup. We take $A=\mathbb{K}\Omega$, with its usual bialgebra structure: the product $m$ obtained by linearization of $\star$ and the coproduct $\Delta$ defined by
$\forall \alpha \in \Omega$, $\Delta(\alpha) =\alpha \otimes \alpha$.
Then $(A,m,\Delta)$ is a counitary bialgebra, unitary if and only if $\Omega$ is a monoid. In~$\leas(A,m,\Delta)$, for any $\alpha,\beta \in \Omega$,
$\Phi(\alpha \otimes \beta)=\alpha\star \beta \otimes \alpha$.
We recover the linearization of $\eas(\Omega,\star)$.
\item[(2)] Let $A$ be a vector space, $1_A\in A$ and $\varepsilon\in A^*$ such that $\varepsilon(1_A)=1$. We define a product and a coproduct on $A$ by
$
\forall a,b\in A$, $ a\cdot b =\varepsilon(a) b$, $
\forall a\in A$, $ \Delta(a) =1_A\otimes a$.
Then $(A,m,\Delta)$ is a bialgebra, with a left unit $1_A$ and a left counit $\varepsilon$. It is unitary if and only if $A$ is one-dimensional; it is counitary if, and only if, $A$ is one-dimensional. In $\leas(A,m,\Delta)$, for any $a,b\in A$, $\Phi(a\otimes b)=b\otimes a$.
\end{enumerate}\end{Example}

\begin{Proposition}\label{prop5.2}
Let $(A,m,\Delta)$ be a bialgebra, not necessarily unitary nor counitary.
\begin{enumerate}\itemsep=0pt
\item[$(1)$] Let us consider the following conditions:
\begin{enumerate}\itemsep=0pt
\item[$(a)$] $\leas(A,m,\Delta)$ is an $\ell$CEDS.
\item[$(b)$] For any $a,b,c\in A$, $\sum \sum a^{(1)}b^{(1)}c\otimes a^{(2)}\otimes b^{(2)}=\sum \sum b^{(1)}a^{(1)}c\otimes a^{(2)}\otimes b^{(2)}$.
\item[$(c)$] For any $a,b,c\in A$, $abc=bac$.
\item[$(d)$] $m$ is commutative.
\end{enumerate}
Then $(d)\Longrightarrow (c)\Longrightarrow (b)\Longleftrightarrow (a)$. If $(A,\Delta)$ has a right counit, then $(c)\Longleftrightarrow (a)$. If $(A,m,\Delta)$ has a right counit and a right unit, then $(d)\Longleftrightarrow (a)$.
\item[$(2)$] Let us consider the following conditions:
\begin{enumerate}\itemsep=0pt
\item[$(a)$] $\leas(A,m,\Delta)$ is a dual $\ell$CEDS.
\item[$(b)$] For any $a,b,c\in A$,
$\sum a^{(1)}b\otimes a^{(2)}c\otimes a^{(3)}=\sum a^{(2)}b\otimes a^{(1)}c \otimes a^{(3)}$.
\item[$(c)$] $(\Delta \otimes \operatorname{Id})\circ \Delta=(\tau \otimes \operatorname{Id})\circ (\Delta \otimes \operatorname{Id})\circ \Delta$.
\item[$(d)$] $\Delta$ is cocommutative.
\end{enumerate}
Then $(d)\Longrightarrow (c)\Longrightarrow (b)\Longleftrightarrow (a)$. If $(A,m)$ has a right unit, then $(c)\Longleftrightarrow (a)$. If $(A,m,\Delta)$ has a right counit and a right unit, then $(d)\Longleftrightarrow (a)$.
\end{enumerate}
\end{Proposition}

\begin{proof}
(1) Obviously, $(d)\Longrightarrow (c)\Longrightarrow (b)$. Let $a,b,c\in A$. Then
\begin{align*}
(\operatorname{Id} \otimes \Phi)\circ (\operatorname{Id} \otimes \tau)\circ (\tau \otimes \operatorname{Id})\circ (\Phi\otimes \operatorname{Id})(b\otimes c\otimes a)&=\sum \sum b^{(2)}\otimes a^{(1)}b^{(1)}c\otimes a^{(2)},\\
(\tau \otimes \operatorname{Id})\circ (\Phi\otimes \operatorname{Id})\circ (\operatorname{Id} \otimes \Phi)\circ (\operatorname{Id} \otimes \tau)(b\otimes c\otimes a)&=\sum \sum b^{(2)}\otimes b^{(1)}a^{(1)}c\otimes a^{(2)},
\end{align*}
so $(a)\Longleftrightarrow (b)$. If $(b)$ is satisfied and if $(A,\Delta)$ has a right counit $\varepsilon$, applying $(\operatorname{Id} \otimes \varepsilon \otimes \varepsilon)$ to $(b)$, we obtain $(c)$.
If $(c)$ is satisfied and $(A,m)$ has a right unit $1_A$, taking $c=1_A$ in $(c)$, we obtain~$(d)$.

(2) Obviously, $(d)\Longrightarrow (c)\Longrightarrow (b)$. Let $a,b,c\in A$. Then
\begin{align*}
(\Phi\otimes \operatorname{Id})\circ (\tau \otimes \operatorname{Id})\circ (\operatorname{Id} \otimes \tau)\circ (\operatorname{Id} \otimes \Phi)(b\otimes a\otimes c)&=\sum a^{(2)}b\otimes a^{(3)}\otimes a^{(1)}c,\\
(\operatorname{Id} \otimes \tau)\circ (\operatorname{Id} \otimes \Phi)\circ (\Phi\otimes \operatorname{Id})\circ (\tau \otimes \operatorname{Id})(b\otimes a\otimes c)&=\sum a^{(1)}b\otimes a^{(3)}\otimes a^{(2)}c,
\end{align*}
so $(a)\Longleftrightarrow (b)$. If $(b)$ is satisfied and if $(A,m)$ has a right unit $1_A$, taking $b=c=1_A$ in $(b)$, we~obtain~$(c)$. If $(c)$ is satisfied
and $(A,\Delta)$ has a~right counit $\varepsilon$, applying $(\operatorname{Id} \otimes \operatorname{Id} \otimes \varepsilon)$ to $(c)$, we~obtain~$(d)$. \end{proof}

\begin{Proposition}
Let $(A,m,\Delta)$ be a finite-dimensional bialgebra, not necessarily unitary nor counitary. Then $\leas(A,m,\Delta)^*=\leas(A^*,\Delta^*,m^*)$.
\end{Proposition}

\begin{proof}
Let $f,g\in A^*$. For any $a,b\in A$,
\begin{align*}
\Phi^*&=((m\otimes \operatorname{Id}_A)\circ (\operatorname{Id}_A\otimes \tau)\circ (\Delta\otimes \operatorname{Id}_A) )^*=(\Delta^*\otimes \operatorname{Id}_{A^*})\circ (\operatorname{Id}_{A^*}\otimes \tau)\circ (m^*\otimes \operatorname{Id}_{A^*}).
\end{align*}
Therefore, $\leas(A,m,\Delta)^*=\leas(A^*,\Delta^*,m^*)$. \end{proof}

\begin{Proposition}\label{prop5.4}
Let $(A,m,\Delta)$ be a bialgebra.
\begin{enumerate}\itemsep=0pt
\item[$(1)$] We assume that $(A,m)$ has a right unit $1_A$.
\begin{itemize}\itemsep=0pt
\item If $1_A$ is not a unit of $(A,m)$, the unique left unit of $\leas(A,m,\Delta)$ is $0$. If $1_A$ is a unit of $(A,m)$, then the left units of $\leas(A,m,\Delta)$ are the elements $a\in A$ such that $\Delta(a)=1_A\otimes a$.
\end{itemize}
\item[$(2)$] We assume that $(A,\Delta)$ has a right counit $\varepsilon_A$.
\begin{itemize}\itemsep=0pt
\item If $\varepsilon_A$ is not a unit of $(A,\Delta)$, the unique left counit of $\leas(A,m,\Delta)$ is $0$. If $\varepsilon_A$ is a counit of $(A,\Delta)$, then the left counits of $\leas(A,m,\Delta)$ are the elements $\lambda \in A^*$
such that $\lambda \circ m=\varepsilon\otimes \lambda$.
\end{itemize}\end{enumerate}\end{Proposition}

\begin{proof}
(1) Let us assume that $\leas(A,m,\Delta)$ has a nonzero left unit $a$. Let us choose $\lambda \in A^*$ such that $\lambda(a)=1$. For any $b\in A$,
\begin{align*}
(\operatorname{Id} \otimes \lambda)\circ \Phi(a\otimes b)&=\underbrace{\left(\sum a^{(1)}\lambda\bigl(a^{(2)}\bigr)\right)}_{=a'} b=(\operatorname{Id} \otimes \lambda)(b\otimes a)=b\lambda(a)=b,
\end{align*}
so $a'$ is a left unit of $(A,m)$. Then $a'1_A=a'=1_A$, so $a'=1_A$ is a unit. Moreover, for $b=1_A$,
\[\Phi(a\otimes 1_A)=\sum a^{(1)}1_A\otimes a^{(2)}=\Delta(a)=1_A\otimes a.\]
Conversely, if $1_A$ is a unit of $(A,m)$ and $\Delta(a)=1_A\otimes a$, then $a$ is clearly a left unit of $\leas(A,m,\Delta)$.

(2) Let us assume that $\leas(A,m,\Delta)$ has a nonzero left counit $\lambda$. Let us choose $b\in A$ such that $\lambda(b)=1$. For any $a\in A$,
\begin{align*}
(\lambda \otimes \operatorname{Id})\circ \Phi(a\otimes b)&=\sum \lambda\bigl(a^{(1)}b\bigr)a^{(2)}=a\lambda(b)=a.
\end{align*}
If we define $\lambda'\colon A\longrightarrow \mathbb{K}$ by $\lambda'(a)=\lambda(ab)$, then $\lambda'$ is a left counit of $(A,\Delta)$. As $\varepsilon_A$ is a right counit of $(A,\Delta)$,
$(\lambda'\otimes \varepsilon_A)\circ \Delta=\lambda'=\varepsilon_A$,
so $\lambda'=\varepsilon_A$ is a counit of $(A,\Delta)$. Moreover, for any~${a,b\in A}$,
\begin{align*}
(\lambda \otimes \varepsilon_A)\circ \Phi(a\otimes b)&=\varepsilon_A(a)\lambda(b)=\sum \lambda\bigl(a^{(1)}b\bigr)\varepsilon_A\bigl(a^{(2)}\bigr)=\lambda(ab),
\end{align*}
so $\lambda \circ m=\varepsilon_A\otimes \lambda$. Conversely, if $\varepsilon_A$ is a counit of $(A,\Delta)$ and $\lambda \circ m=\varepsilon_A\otimes \lambda$, then for any $a,b\in A$,
\begin{align*}
(\lambda \otimes \operatorname{Id})\circ \Phi(a\otimes b)&=\sum \lambda\bigl(a^{(2)}b\bigr)a^{(2)}=\sum \varepsilon_A\bigl(a^{(1)}\bigr)\lambda(b)a=\lambda(b)a,
\end{align*}
so $\lambda$ is a left counit of $\leas(A,m,\Delta)$. \end{proof}

More generally, we can obtain other $\ell$EAS with the help of a bialgebra projection or with certain linear forms.

\begin{Proposition}
Let $(A,m,\Delta)$ be a bialgebra, not necessarily unitary nor counitary, and $\pi\colon A\longrightarrow A$ be a bialgebra morphism such that $\pi^2=\pi$. For any $a,b\in A$, we define $\Phi\colon A\otimes A\longrightarrow A\otimes A$~by
\begin{align*}
\Phi(a\otimes b)=(m\otimes \pi)\circ (\operatorname{Id}_A\otimes \tau)\circ (\Delta\otimes \operatorname{Id}_A)(a\otimes b)=\sum a^{(1)}b\otimes \pi\bigl(a^{(2)}\bigr).
\end{align*}
 Then $(A,\Phi)$ is an $\ell$EAS.
\end{Proposition}

\begin{proof}
We define $\delta=(\operatorname{Id} \otimes \pi)\circ \Delta$. Then $(A,m,\delta)$ is a bialgebra. Note that it is not counitary, except if $(A,\Delta)$ is counitary and $\pi=\operatorname{Id}_A$. We can then apply Proposition \ref{prop5.1} to $(A,m,\delta)$.
\end{proof}

\begin{Example}
Let $(\Omega,\star)$ be a semigroup and $\pi\colon\Omega\longrightarrow \Omega$ be a semigroup morphism such that~${\pi^2=\pi}$. We take $A=\mathbb{K}\Omega$, with its usual bialgebra structure. Then in $\leas(A,m,\Delta)$, for any $\alpha,\beta \in \Omega$,
$\Phi(\alpha \otimes \beta)=\alpha\star \beta \otimes \pi(\alpha)$.
We recover the linearization of $\eas(\Omega,\star,\pi)$.
\end{Example}

\begin{Proposition}
Let $(A,\Delta)$ be a coalgebra, not necessarily counitary, and $f\in A^*$ such that $(f\otimes f)\circ \Delta=f$. We put, for any $a,b\in A$,
\[\Phi(a\otimes b)=\sum f\bigl(a^{(1)}\bigr)b\otimes a^{(2)}.\]
Then $(A,\Phi)$ is an $\ell$CEDS.
\end{Proposition}

\begin{proof}
We define a product on $A$ by $a\star b=f(a)b$. It is associative. Moreover, for any $a,b\in A$, as $(f\otimes f)\circ \Delta=f$,
\begin{align*}
\Delta(a\star b)&=f(a)\sum b^{(1)}\otimes b^{(2)}=\sum f\bigl(a^{(1)}\bigr)f\bigl(a^{(2)}\bigr)\sum b^{(1)}\otimes b^{(2)}\\
&=\sum\sum f\bigl(a^{(1)}\bigr)b^{(1)}\otimes f\bigl(a^{(2)}\bigr)b^{(2)}=\Delta(a)\star \Delta(b),
\end{align*}
so $(A,\star,\Delta)$ is a bialgebra, and $(A,\Phi)=\leas(A,\star,\Delta)$. Moreover, for any $a,b,c\in A$,
\[a\star b\star c=f(a)f(b)c=f(b)f(c)a=b\star a\star c.\]
By Proposition \ref{prop5.2}, $(A,\Phi)$ is an $\ell$CEDS. \end{proof}

\begin{Example}
Let $\Omega$ be a set, $A=\mathbb{K}\Omega$ be the associated coalgebra (where any $\alpha \in \Omega$ is a~group-like element), and $\Omega'\subseteq \Omega$ be any set. We define the linear form $f\colon A\longrightarrow \mathbb{K}$ by
\begin{align*}
 \forall \alpha \in \Omega, f(\alpha) =\begin{cases}
1& \text{if }\alpha \in \Omega',\\
0 & \text{otherwise}.
\end{cases} \end{align*}
For any $\alpha \in \Omega$, $(f\otimes f)\circ \Delta(\alpha)=f(\alpha)^2=f(\alpha)$, so we obtain an $\ell$CEDS such that for any $\alpha,\beta \in \Omega$,
\[\Phi(\alpha \otimes \beta)=\begin{cases}
\beta\otimes \alpha&\text{if }\alpha \in \Omega',\\
0 & \text{otherwise}.
\end{cases}\] \end{Example}

\subsection[A functor from Hopf algebras to lEAS]{A functor from Hopf algebras to $\boldsymbol{\ell}$EAS}

\begin{Proposition}\label{prop5.7}
Let $(A,m,\Delta)$ be a Hopf algebra, of antipode $S$. For any $a,b\in A$, we define $\Phi\colon A\otimes A\longrightarrow A\otimes A$ by
\begin{gather*}
\Phi(a\otimes b)=(\operatorname{Id}_A\otimes m)\circ (\operatorname{Id}_A\otimes S\otimes \operatorname{Id}_A)\circ (\Delta \otimes \operatorname{Id})\circ \tau(a\otimes b)=\sum b^{(1)}\otimes S\bigl(b^{(2)}\bigr)a.
\end{gather*}
Then $(A,\Phi)$ is an $\ell$EAS, denoted by $\leas'(A,m,\Delta)$. It is non-degenerate, and $\bigl(A,\Phi^{-1}\bigr)=\leas(A,m,\Delta^{\rm op})$.
\end{Proposition}

\begin{proof}
Let $a,b,c\in A$. Then
\begin{gather*}
(\operatorname{Id} \otimes \Phi)\circ(\Phi\otimes \operatorname{Id})\circ (\operatorname{Id} \otimes \Phi)(a\otimes b\otimes c)\\
\qquad=\sum \sum c^{(1)}\otimes S\bigl(c^{(3)}\bigr)^{(1)}b^{(1)}\otimes S\bigl(S\bigl(c^{(3)}\bigr)^{(2)}b^{(2)}\bigr)S\bigl(c^{(2)}\bigr)a\\
\qquad=\sum \sum c^{(1)}\otimes S\bigl(c^{(4)}\bigr)b^{(1)}\otimes S\bigl(S\bigl(c^{(3)}\bigr)b^{(2)}\bigr)S\bigl(c^{(2)}\bigr)a\\
\qquad=\sum \sum c^{(1)}\otimes S\bigl(c^{(4)}\bigr)b^{(1)}\otimes S\bigl(c^{(2)}S\bigl(c^{(3)}\bigr)b^{(2)}\bigr)a\\
\qquad=\sum \sum c^{(1)}\otimes S\bigl(c^{(2)}\bigr)b^{(1)}\otimes S\bigl(b^{(2)}\bigr)a\\
\qquad=(\Phi\otimes \operatorname{Id})\circ (\operatorname{Id} \otimes \tau)\circ (\Phi\otimes \operatorname{Id})(a\otimes b\otimes c),
\end{gather*}
so $(A,\Phi)$ is an $\ell$EAS.

Let $(A,\Psi)=\leas(A,m,\Delta^{\rm op})$: for any $a,b\in A$, $\Psi(a\otimes b)=\sum a^{(2)}b\otimes a^{(1)}$. Then
\begin{gather*}
\Phi\circ \Psi(a \otimes b)=\sum a^{(1)}\otimes S\bigl(a^{(2)}\bigr)a^{(3)}b=a\otimes b,\\
\Psi\circ \Phi(a\otimes b)=\sum b^{(1)}S\bigl(b^{(2)}\bigr)a\otimes b^{(3)}=a\otimes b,
\end{gather*}
so $\Phi$ is bijective, of inverse $\Psi$. \end{proof}

\begin{Example}
Let $(G,\star)$ be a group and let $A=\mathbb{K} G^{\rm op}$ be the Hopf algebra of the opposite of this group. A basis of $\leas'(A,m,\Delta)$ is given by $G$ itself and, for any $\alpha,\beta\in G$,
$\Phi(\alpha \otimes \beta)=\beta \otimes \alpha \star \beta^{-1}$.
We recover in this way the linearization of $\eas'(G,\star)$.
\end{Example}

\begin{Corollary}
Let $(A,m,\Delta)$ be a bialgebra, such that $(A,m,\Delta^{\rm op})$ is a bialgebra. Then $(A,\Phi)=\leas(A,m,\Delta)$ is non-degenerate and $\bigl(A,\Phi^{-1}\bigr)=\leas'(A,m,\Delta^{\rm op})$.
\end{Corollary}

\begin{Proposition}
Let $(A,m,\Delta)$ be a Hopf algebra.
\begin{enumerate}\itemsep=0pt
\item[$(1)$] Then $\leas'(A,m,\Delta)$ is an $\ell$CEDS if and only if $\Delta \circ S=\Delta^{\rm op}\circ S$.
\item[$(2)$] Then $\leas'(A,m,\Delta)$ is a dual $\ell$CEDS if and only if $S\circ m=S\circ m^{\rm op}$.
\end{enumerate}
\end{Proposition}

\begin{proof}
(1) Let $a,b,c\in A$.
\begin{align}
\label{eq34} (\operatorname{Id} \otimes \Phi)\circ (\operatorname{Id} \otimes \tau)\circ (\tau \otimes \operatorname{Id})\circ (\Phi\otimes \operatorname{Id})(b\otimes c\otimes a)&=\sum S\bigl(b^{(3)}\bigr)a\otimes b^{(1)}\otimes S\bigl(b^{(2)}\bigr)c,\\
\nonumber (\tau \otimes \operatorname{Id})\circ (\Phi\otimes \operatorname{Id})\circ (\operatorname{Id} \otimes \Phi)\circ (\operatorname{Id} \otimes \tau)(b\otimes c\otimes a)&=\sum S\bigl(b^{(2)}\bigr)a\otimes b^{(1)}\otimes S\bigl(b^{(3)}\bigr)c.
\end{align}
If $\Delta \circ S=\Delta^{\rm op}\circ S$, then
\begin{align*}
\sum b^{(1)}\otimes S\bigl(b^{(2)}\bigr)\otimes S\bigl(b^{(3)}\bigr)&=\sum \sum b^{(1)}\otimes S\bigl(b^{(2)}\bigr)^{(2)} \otimes S\bigl(b^{(2)}\bigr)^{(1)}\\
&=\sum \sum b^{(1)}\otimes S\bigl(b^{(2)}\bigr)^{(1)} \otimes S\bigl(b^{(2)}\bigr)^{(2)}\\
&=\sum b^{(1)}\otimes S\bigl(b^{(3)}\bigr)\otimes S\bigl(b^{(2)}\bigr),
\end{align*}
which implies that $(A,\Phi)$ is an $\ell$CEDS. Conversely, taking $a=c=1_A$, we obtain, in \eqref{eq34},
\[\sum S\bigl(b^{(3)}\bigr)\otimes b^{(1)}\otimes S\bigl(b^{(2)}\bigr)=\sum S\bigl(b^{(2)}\bigr)\otimes b^{(1)}\otimes S\bigl(b^{(3)}\bigr).\]
Applying $\operatorname{Id} \otimes \varepsilon \otimes \operatorname{Id}$, we obtain
\begin{align*}
\Delta \circ S(b)&=\sum S(b)^{(1)}\otimes S(b)^{(2)}=\sum S\bigl(b^{(2)}\bigr)\otimes S\bigl(b^{(1)}\bigr)\\
&=\sum S\bigl(b^{(1)}\bigr)\otimes S\bigl(b^{(2)}\bigr)=\sum S(b)^{(2)}\otimes S(b)^{(1)}=\Delta^{\rm op}\circ S(b).
\end{align*}

(2) Let $a,b,c\in A$. Then
\begin{gather*}
(\Phi\otimes \operatorname{Id})\circ (\tau \otimes \operatorname{Id})\circ (\operatorname{Id} \otimes \tau)\circ (\operatorname{Id} \otimes \Phi)(b\otimes a\otimes c)=\sum \sum a^{(1)}\otimes S\bigl(a^{(2)}\bigr)S\bigl(c^{(2)}\bigr)b\otimes c^{(1)},\\
\nonumber(\operatorname{Id} \otimes \tau)\circ (\operatorname{Id} \otimes \Phi)\circ (\Phi\otimes \operatorname{Id})\circ (\tau \otimes \operatorname{Id})(b\otimes a\otimes c)=\sum \sum a^{(1)}\otimes S\bigl(c^{(2)}\bigr)S\bigl(a^{(2)}\bigr)b\otimes c^{(1)}.
\end{gather*}
If $S\circ m=S\circ m^{\rm op}$, then $m^{\rm op}\circ (S\otimes S)=m\circ (S\otimes S)$, which implies that $(A,\Phi)$ is a dual $\ell$CEDS. Conversely, taking $b=1$ and applying $\varepsilon \otimes \operatorname{Id} \otimes \varepsilon$, we obtain
$S(a)S(c)=S(c)S(a)$,
so~${m^{\rm op}\circ S=m\circ (S\otimes S)=m^{\rm op}\circ (S\otimes S)=S\circ m}$. \end{proof}

\begin{Remark}
In particular, if $S$ is invertible, then $\leas'(A,m,\Delta)$ is an $\ell$CEDS if and only if $(A,m)$ is commutative; it is a dual $\ell$CEDS if and only if $(A,\Delta)$ is cocommutative.
\end{Remark}

\begin{Proposition}
Let $(A,m,\Delta)$ be a finite-dimensional Hopf algebra. Then \[
\leas'(A,m,\Delta)^* =\leas'(A^*,\Delta^{*,{\rm op}},m^{*,{\rm op}}).\]
\end{Proposition}

\begin{proof}
Let $f,g\in A^*$. For any $a,b\in A$,
\begin{align*}
\Phi^*(f\otimes g)(a\otimes b)&=(f\otimes g)(\Phi(a\otimes b))=\sum (f\otimes g)\bigl(b^{(1)}\otimes S\bigl(b^{(2)}a\bigr)\bigr)\\
&=\sum \sum \bigl(f\otimes g^{(1)}\otimes g^{(2)}\bigr)\bigl(b^{(1)}\otimes S\bigl(b^{(2)}\bigr)\otimes a\bigr)\\
&=\sum \bigl(f\otimes S^*\bigl(g^{(1)}\bigr)\otimes g^{(2)}\bigr)\bigl(b^{(1)}\otimes b^{(2)}\otimes a\bigr)\\
&=\sum \bigl(g^{(2)}\otimes fS^*\bigl(g^{(1)}\bigr)\bigr)(a\otimes b),
\end{align*}
so $\Phi^*(f\otimes g)=\sum g^{(2)}\otimes fS^*\bigl(g^{(1)}\bigr)$, which is the $\ell$EAS attached to the Hopf algebra $(A^*,\Delta^{*,{\rm op}},\allowbreak m^{*,{\rm op}})$, whose antipode is $S^*$. \end{proof}

Recall from \cite{Sweedler1969-2} that a right integral of a Hopf algebra $(A,m,\Delta)$ is a linear map $f\in A^*$ such that for any $\mu\in A^*$,
$(\lambda \otimes \mu)\circ \Delta=\mu(1_A) \lambda$.

\begin{Proposition} \label{prop5.11}
Let $(A,m,\Delta)$ be a Hopf algebra.
\begin{enumerate}\itemsep=0pt
\item[$(1)$] Let $a\in A$. It is a left unit of $\leas'(A,m,\Delta)$ if and only if for any $b\in A$, $S(b)a=\varepsilon(b)a$.
\item[$(2)$] Let $\lambda \in A^*$. It is a left counit of $\leas'(A,m,\Delta)$ if and only if for any $a\in A$, \linebreak $\sum \lambda\bigl(b^{(1)}\bigr) S\bigl(b^{(2)}\bigr)=\lambda(b)1_A$.
In particular, right integrals on $(A,m,\Delta)$ are left counit of $\leas'(A,m,\Delta)$; if $S$ is invertible, then the converse is true.
\end{enumerate}\end{Proposition}

\begin{proof}
(1) Let $a\in A$. Then its a left unit if and only if for any $b\in A$, $\sum b^{(1)}\otimes S\bigl(b^{(2)}\bigr)a=b\otimes a$. Applying $\varepsilon\otimes \operatorname{Id}$, if $a$ is a left unit, for any $b\in B$, $S(b)a=\varepsilon(b)\otimes a$.
Conversely, if this holds, then for any $b\in B$,
\[\Phi(a\otimes b)=\sum b^{(1)}\otimes S\bigl(b^{(1)}\bigr)a=\sum b^{(1)}\otimes \varepsilon\bigl(b^{(2)}\bigr) a=b\otimes a.\]

(2) Let $\lambda \in A^*$. It is a left counit if and only if for any $a,b\in A$,
\[\sum \lambda\bigl(b^{(1)}\bigr)S\bigl(b^{(2)}\bigr)a=a\lambda(b).\]
If $\lambda$ is a left counit, taking $a=1_A$, we obtain that for any $b\in A$,
$\lambda\bigl(b^{(1)}\bigr)S\bigl(b^{(2)}\bigr)=\lambda(b)1_A$. Conversely, if this holds, then for any $a,b\in A$,{\samepage
\[(\lambda\otimes \operatorname{Id})\circ \Phi(a\otimes b)=\sum \lambda\bigl(b^{(1)}\bigr)S\bigl(b^{(2)}\bigr)a=\lambda(b)a=(\operatorname{Id} \otimes \lambda)(a\otimes b),\]
so $\lambda$ is a left counit.}

Let us assume that $\lambda$ is a right integral of $(A,m,\Delta)$. For any $b\in A$, for any $\mu\in A^*$,
\begin{align*}
\sum \lambda\bigl(b^{(1)}\bigr)\mu\bigl(S\bigl(b^{(2)}\bigr)\bigr)&=(\lambda\otimes \mu\circ S)\circ \Delta(b)=\mu\circ S(1_A)\lambda(b)=\mu(1_A)\lambda(b).
\end{align*}
As this holds for any $\mu\in A^*$, $\sum \lambda\bigl(b^{(1)}\bigr)S\bigl(b^{(2)}\bigr)=\lambda(b)1_A$, so $\lambda$ is a right integral. Let us now assume that $S$ is invertible and that $\lambda$ is a left counit.
Let $\nu \in A^*$. For any $b\in A$, if $\mu=\nu \circ S^{-1}$,
\begin{align*}
\sum \lambda\bigl(b^{(1)}\bigr)\nu\bigl(b^{(2)}\bigr)&=\sum \lambda\bigl(b^{(1)}\bigr)\mu\circ S\bigl(b^{(2)}\bigr)=\lambda(b)\mu(1_A)\\
&=\lambda(b)\nu\circ S^{-1}(1_A)=\lambda(b)\nu(1_A).
\end{align*}
So $\lambda$ is a right integral.
\end{proof}

\subsection{From left units and counits to bialgebras}

\begin{Theorem}\label{theo5.12}
Let $(A,\Phi)$ be an $\ell$EAS.
\begin{enumerate}\itemsep=0pt
\item[$(1)$] If $a$ is a special vector of eigenvalue $1$ of $(A,\Phi)$, then $\Delta_a\colon A\longrightarrow A\otimes A$ defined by $\Delta_a(b)=\Phi(b\otimes a)$ is a coassociative coproduct.
\item[$(2)$] If $\varepsilon$ is a special vector of eigenvalue $1$ of $(A,\Phi)^*$, that is to say if $(\varepsilon\otimes \varepsilon)\circ \Phi=\varepsilon\otimes \varepsilon$,
 then $m_\varepsilon\colon A\otimes A\longrightarrow A$ defined by $m_\varepsilon=(\operatorname{Id} \otimes \varepsilon)\circ \Phi$ is an associative product.
\item[$(3)$] If $a$ is a left unit of $(A,\Phi)$ and $\varepsilon$ is a left counit of $(A,\Phi)$ such that $\varepsilon(a)=1$, then $(A,m_\varepsilon,\Delta_a)$ is a bialgebra, with $a$ as a left unit and $\varepsilon$ as a left counit.
Moreover, $(A,\Phi)=\leas(A,m_\varepsilon,\Delta_a)$.
\end{enumerate}\end{Theorem}

\begin{proof}
(1) For any $b\in A$,
\begin{gather*}
(\operatorname{Id} \otimes \Phi)\circ (\Phi\otimes \operatorname{Id})\circ (\operatorname{Id} \otimes \Phi)(b\otimes a\otimes a)=(\operatorname{Id} \otimes \Phi)\circ (\Phi\otimes \operatorname{Id})(b\otimes a\otimes a)\\
\phantom{(\operatorname{Id} \otimes \Phi)\circ (\Phi\otimes \operatorname{Id})\circ (\operatorname{Id} \otimes \Phi)(b\otimes a\otimes a)}{}=(\operatorname{Id} \otimes \Phi)(\Delta_a(b)\otimes a)=(\operatorname{Id} \otimes \Delta_a)\circ \Delta_a(b),\\
(\Phi\otimes \operatorname{Id})\circ (\operatorname{Id} \otimes \tau)\circ (\Phi\otimes \operatorname{Id})(b\otimes a\otimes a)=(\Phi\otimes \operatorname{Id})\circ (\operatorname{Id} \otimes \tau)(\Delta_a(b)\otimes a)\\
\phantom{(\Phi\otimes \operatorname{Id})\circ (\operatorname{Id} \otimes \tau)\circ (\Phi\otimes \operatorname{Id})(b\otimes a\otimes a)}{}=(\Delta_a\otimes \operatorname{Id})\circ \Delta_a(b).
\end{gather*}
Hence, $\Delta_a$ is coassociative.

(2) We obtain, as $\varepsilon$ is a special vector of eigenvalue $1$ of $(A,\Phi)^*$,
\begin{gather*}
(\operatorname{Id} \otimes \varepsilon\otimes \varepsilon)\circ (\operatorname{Id} \otimes \Phi)\circ (\Phi\otimes \operatorname{Id})\circ (\operatorname{Id} \otimes \Phi)=(\operatorname{Id} \otimes \varepsilon\otimes \varepsilon)\circ (\Phi\otimes \operatorname{Id})\circ (\operatorname{Id} \otimes \Phi)\\
\phantom{(\operatorname{Id} \otimes \varepsilon\otimes \varepsilon)\circ (\operatorname{Id} \otimes \Phi)\circ (\Phi\otimes \operatorname{Id})\circ (\operatorname{Id} \otimes \Phi)}{}=(\operatorname{Id} \otimes \varepsilon)\circ \Phi\circ ((\operatorname{Id} \otimes \varepsilon)\circ \Phi))\\
\phantom{(\operatorname{Id} \otimes \varepsilon\otimes \varepsilon)\circ (\operatorname{Id} \otimes \Phi)\circ (\Phi\otimes \operatorname{Id})\circ (\operatorname{Id} \otimes \Phi)}{}=m_\varepsilon\circ(\operatorname{Id} \otimes m_\varepsilon),\\
(\operatorname{Id} \otimes \varepsilon\otimes \varepsilon)\circ (\Phi\otimes \operatorname{Id})\circ (\operatorname{Id} \otimes \tau)\circ (\Phi\otimes \operatorname{Id})=(\operatorname{Id} \otimes \varepsilon )\circ \Phi\circ ((\operatorname{Id} \otimes \varepsilon )\otimes \operatorname{Id})\\
\phantom{(\operatorname{Id} \otimes \varepsilon\otimes \varepsilon)\circ (\Phi\otimes \operatorname{Id})\circ (\operatorname{Id} \otimes \tau)\circ (\Phi\otimes \operatorname{Id})}{}=m_\varepsilon\circ(m_\varepsilon\otimes \operatorname{Id}).
\end{gather*}
As a consequence, $m_\varepsilon$ is associative.

(3) As $a$ is a left unit, it is a special vector of eigenvalue $1$ of $(A,\Phi)$, so $\Delta_a$ is coassociative. Moreover, for any $b\in A$,
\[(\varepsilon\otimes \operatorname{Id})\circ \Delta_a(b)=(\varepsilon\otimes \operatorname{Id})\circ \Phi(b\otimes a)=(\operatorname{Id} \otimes \varepsilon)(b\otimes a)=b\varepsilon(a)=b,\]
so $\varepsilon$ is a left counit of $\Delta_a$. As $\varepsilon$ is a left counit, it is a special vector of eigenvalue $1$ of $(A,\Phi)^*$, so $m_\varepsilon$ is associative. Moreover, for any $b\in A$,{\samepage
\[m_\varepsilon(a\otimes b)=(\operatorname{Id} \otimes \varepsilon)\circ \Phi(a\otimes b)=(\operatorname{Id} \otimes \varepsilon)(b\otimes a)=b\varepsilon(a)=b,\]
so $a$ is a left unit of $m_\varepsilon$.}

Let $b_1,b_2\in A$.
\begin{align*}
\Delta_a(b_1b_2)={}&(\operatorname{Id} \otimes \varepsilon\otimes \operatorname{Id} \otimes \varepsilon)\circ (\Phi\otimes \Phi)\circ (\operatorname{Id} \otimes \tau\otimes \operatorname{Id})\circ (\Phi\otimes \Phi)\\
&\circ (\operatorname{Id} \otimes \tau\otimes \operatorname{Id})(b_1\otimes b_2\otimes a\otimes a)\\
={}&(\operatorname{Id} \otimes \varepsilon\otimes \operatorname{Id} \otimes \varepsilon)\circ (\operatorname{Id} \otimes \operatorname{Id} \otimes \Phi)\circ (\Phi\otimes \operatorname{Id} \otimes \operatorname{Id})\circ (\operatorname{Id} \otimes \tau\otimes \operatorname{Id})\circ (\Phi\otimes \operatorname{Id} \otimes \operatorname{Id})\\
&\circ(\operatorname{Id} \otimes \operatorname{Id} \otimes \Phi)\circ (\operatorname{Id} \otimes \tau\otimes \operatorname{Id})(b_1\otimes b_2\otimes a\otimes a)\\
={}&(\operatorname{Id} \otimes \varepsilon\otimes \operatorname{Id} \otimes \varepsilon)\circ (\operatorname{Id} \otimes \operatorname{Id} \otimes \Phi)\circ (\operatorname{Id}\otimes \Phi \otimes \operatorname{Id})\circ (\Phi \otimes \operatorname{Id}\otimes \operatorname{Id})\circ (\operatorname{Id}\otimes \Phi\otimes \operatorname{Id})\\
&\circ(\operatorname{Id} \otimes \operatorname{Id} \otimes \Phi)\circ (\operatorname{Id} \otimes \tau\otimes \operatorname{Id})(b_1\otimes b_2\otimes a\otimes a)\\
={}&(\operatorname{Id} \otimes \operatorname{Id} \otimes \varepsilon)\circ (\operatorname{Id} \otimes \Phi)\circ (\operatorname{Id}\otimes ((\varepsilon \otimes \operatorname{Id})\circ\Phi) \otimes \operatorname{Id})\circ (\Phi \otimes \operatorname{Id}\otimes \operatorname{Id})\\
&\circ (\operatorname{Id}\otimes \Phi\otimes \operatorname{Id})\circ(\operatorname{Id} \otimes \operatorname{Id} \otimes \Phi)\circ (\operatorname{Id} \otimes \tau\otimes \operatorname{Id})(b_1\otimes b_2\otimes a\otimes a)\\
={}&(\operatorname{Id} \otimes \operatorname{Id} \otimes \varepsilon)\circ (\operatorname{Id} \otimes \Phi)\circ (\operatorname{Id}\otimes \operatorname{Id} \otimes \varepsilon \otimes \operatorname{Id})\circ (\Phi \otimes \operatorname{Id}\otimes \operatorname{Id})\\
&\circ (\operatorname{Id}\otimes \Phi\otimes \operatorname{Id})\circ(\operatorname{Id} \otimes \operatorname{Id} \otimes \Phi)\circ (\operatorname{Id} \otimes \tau\otimes \operatorname{Id})(b_1\otimes b_2\otimes a\otimes a)\\
={}&(\operatorname{Id} \otimes \operatorname{Id} \otimes \varepsilon)\circ (\operatorname{Id} \otimes \Phi)\circ (\Phi \otimes \operatorname{Id}\otimes \operatorname{Id})\circ (\operatorname{Id}\otimes (\operatorname{Id}\otimes \varepsilon)\circ\Phi\otimes \operatorname{Id})\\
&\circ(\operatorname{Id} \otimes \operatorname{Id} \otimes \Phi)(b_1\otimes a\otimes b_2\otimes a)\\
={}&(\operatorname{Id} \otimes \operatorname{Id} \otimes \varepsilon)\circ (\operatorname{Id} \otimes \Phi)\circ (\Phi \otimes \operatorname{Id})\circ(\operatorname{Id} \otimes \Phi)(b_1\otimes b_2\otimes a)
\\
={}&\Phi\circ((\operatorname{Id} \otimes \varepsilon)\circ \Phi \otimes \operatorname{Id})(b_1\otimes b_2\otimes a)\\
={}&\Phi(m_\varepsilon(b_1\otimes b_2)\otimes a)
=\Delta_a(m_\varepsilon(b_1\otimes b_2)).
\end{align*}
So $(A,m_\varepsilon,\Delta_a)$ is a bialgebra. Let $(A,\Psi)=\leas(A,m_\varepsilon,\Delta_a)$. For any $b_1,b_2\in A$,
\begin{align*}
\Psi(b_1\otimes b_2)&=(\operatorname{Id} \otimes \varepsilon\otimes \operatorname{Id})\circ (\Phi\otimes \operatorname{Id})\circ (\operatorname{Id} \otimes \tau)\circ (\Phi\otimes \operatorname{Id})(b_1\otimes a\otimes b_2)\\
&=(\operatorname{Id} \otimes \varepsilon\otimes \operatorname{Id})\circ (\operatorname{Id}\otimes \Phi)\circ (\Phi \otimes \operatorname{Id})\circ (\operatorname{Id}\otimes \Phi)(b_1\otimes a\otimes b_2)\\
&=(\operatorname{Id} \otimes \operatorname{Id} \otimes \varepsilon)\circ (\Phi\otimes \operatorname{Id})(b_1\otimes b_2\otimes a)=\Phi(b_1\otimes b_2)\varepsilon(a)=\Phi(b_1\otimes b_2).
\end{align*}
Therefore, $(A,\Phi)=\leas(A,m_\varepsilon,\Delta_a)$. \end{proof}

\begin{Example}
This can be applied for $\ell$EAS $M_{16}$, $M_{17}$ and $M_{18}$ of Example \ref{ex4.1}.
\begin{itemize}
\item For $M_{16}$, taking $a=x$ and $\varepsilon=x^*+y^*$, we obtain
\begin{alignat*}{3}
&\Delta_a(x)=x\otimes x,\qquad&&\Delta_a(y)=y\otimes y,&\\
&m_\varepsilon(x\otimes x)=x,\qquad&& m_\varepsilon(x\otimes y)=y,&\\
&m_\varepsilon(y\otimes x)=y,\qquad&& m_\varepsilon(y\otimes y)=y.&
\end{alignat*}
This is the bialgebra of the semigroup $(\mathbb{Z}/2\mathbb{Z},\times)$, with $x=\overline{1}$ and $y=\overline{0}$:
we recover the linearization of $\mathbf{C3}$.
\item For $M_{17}$, taking $a=x$ and $\varepsilon=x^*+y^*$, we obtain
\begin{alignat*}{3}
&\Delta_a(x)=x\otimes x,\qquad&&\Delta_a(y)=x\otimes x-x\otimes y-y\otimes x+2y\otimes y,&\\
&m_\varepsilon(x\otimes x)=x,\qquad&& m_\varepsilon(x\otimes y)=y,&\\
&m_\varepsilon(y\otimes x)=x,\qquad&& m_\varepsilon(y\otimes y)=y.&
\end{alignat*}
Putting $y'=-x+2y$,\footnote{If the characteristic of the base field $\mathbb{K}$ is not 2.} we obtain
\begin{alignat*}{3}
&\Delta_a(x)=x\otimes x,\qquad&&\Delta_a(y')=y'\otimes y',&\\
&m_\varepsilon(x\otimes x)=x,\qquad&& m_\varepsilon(x\otimes y')=y',&\\
&m_\varepsilon(y'\otimes x)=y',\qquad&& m_\varepsilon(y'\otimes y')=x.&
\end{alignat*}
This is the bialgebra of the semigroup $(\mathbb{Z}/2\mathbb{Z},+)$, with $x=\overline{0}$ and $y=\overline{1}$:
we recover the linearization of $\mathbf{H2}$.
\item For $M_{18}$, we can take any $a\in A$ and any $\varepsilon\in A^*$ such that $\varepsilon(a)=1$. For any $b,c\in A$,
$
\Delta_a(b)=a\otimes b$, $ m_\varepsilon(b\otimes c)=\varepsilon(b)c$.
\end{itemize}\end{Example}

\subsection{Applications to non-degenerate finite CEDS}

From Proposition \ref{prop4.5}, we have the following.

\begin{Proposition}
Let $(\Omega,\rightarrow,\triangleright)$ be a non-degenerate finite CEDS, which we write following Theorem {\rm\ref{theo3.14}} under the form
$(\eas(\Omega_1,*) \rtimes_\succ \eas'(\Omega_2,\star))\times \eas(\Omega_3)$.
Let $g,h\colon\Omega_3\longrightarrow \mathbb{K}$ be two maps such that
\[\sum_{\alpha_3\in \Omega_3} g(\alpha_3)h(\alpha_3)=1.\]
We define a product and a coproduct on $\mathbb{K}\Omega$, putting, for any $(\alpha_1,\alpha_2,\alpha_3)$, $(\beta_1,\beta_2,\beta_3)\in \Omega$,
\begin{gather*}
(\alpha_1,\alpha_2,\alpha_3)\cdot(\beta_1,\beta_2,\beta_3)=\delta_{\alpha_2,\beta_2}g(\alpha_3)(\alpha_1*\beta_1,\beta_2,\beta_3),\\
\Delta(\alpha_1,\alpha_2,\alpha_3)=\sum_{(\beta_2,\beta_3)\in \Omega_2\times \Omega_3}h(\beta_3)(\alpha_1,\beta_2,\beta_3)\otimes \bigl(\beta_2\succ \alpha_1,\alpha_2\star \beta_2^{-1},\alpha_3\bigr).
\end{gather*}
Then $(\mathbb{K}\Omega,\cdot,\Delta)$ is a bialgebra and the linearization of $\Omega$ is $\leas(\mathbb{K}\Omega,\cdot,\Delta)$.
\end{Proposition}

\begin{proof}
By Proposition \ref{prop4.5}, the following is a left unit of $\mathbb{K}\Omega$:
\[a=\sum_{(\alpha_2,\alpha_3)\in \Omega_2\times \Omega_3} h(\alpha_3)(e_1,\alpha_2,\alpha_3),\]
and the following map is a left counit of $\mathbb{K}\Omega$:
\[\varepsilon\colon \ \begin{cases}
\mathbb{K}\Omega \longrightarrow \mathbb{K},\\
(\alpha_1,\alpha_2,\alpha_3) \longrightarrow \delta_{\alpha_2,e_2}g(\alpha_3).
\end{cases}\]
By hypothesis, $\varepsilon(a)=1$. The result comes from a direct application of Theorem~\ref{theo5.12}.
\end{proof}

Similarly, we have the following.

\begin{Proposition}
Let $(\Omega,\rightarrow,\triangleright)$ be a non-degenerate finite dual CEDS, which we write following Corollary {\rm\ref{cor3.15}} under the form
$(\eas(\Omega_2,\star) \ltimes_\prec (\eas'(\Omega_1,*))\times \eas(\Omega_3)$,
Let $g,h\colon\Omega_3\longrightarrow \mathbb{K}$ be two maps such that
\[\sum_{\alpha_3\in \Omega_3} g(\alpha_3)h(\alpha_3)=1.\]
We define a product and a coproduct on $\mathbb{K}\Omega$, putting, for any $(\alpha_1,\alpha_2,\alpha_3)$, $(\beta_1,\beta_2,\beta_3)\in \Omega$,
\begin{gather*}
(\alpha_1,\alpha_2,\alpha_3)\cdot(\beta_1,\beta_2,\beta_3)=\delta_{\alpha_1,\beta_1}g(\alpha_3)(\alpha_2*\beta_2,\beta_1\prec \alpha_2,\beta_3),\\
\Delta(\alpha_2,\alpha_1,\alpha_3)=\sum_{(\beta_1,\beta_3)\in \Omega_1\times \Omega_3}h(\beta_3)(\alpha_2,\beta_1\prec \alpha_2,\beta_3)\otimes \bigl(\alpha_2,\alpha_1\star\bigl(\beta_1^{-1}\prec \alpha_2^{-1}\bigr),\alpha_3\bigr).
\end{gather*}
Then $(\mathbb{K}\Omega,\cdot,\Delta)$ is a bialgebra and the linearization of $\Omega$ is $\leas(\mathbb{K}\Omega,\cdot,\Delta)$.
\end{Proposition}

\subsection{Applications to Hopf algebras of groups}

In all this paragraph, $G$ is a group. We denote by $\mathbb{K} G$ the associated Hopf algebra. If $G$ is finite, we denote by $\mathbb{K}^G$ the Hopf algebra of functions over $G$, with its basis $(\delta_g)_{g\in G}$, dual of the basis~$G$ of $\mathbb{K} G$.

\begin{Corollary}\label{cor5.15}
If $G$ is finite, then $\leas'(\mathbb{K} G)$ is isomorphic to $\leas\bigl(\mathbb{K}^G\bigr)$, and $\leas'\bigl(\mathbb{K}^G\bigr)$ is isomorphic to $\leas(\mathbb{K} G^{\rm op})$.
\end{Corollary}

\begin{proof}
As $G$ is finite, $a= \sum_{g\in G} g$ is a right integral of $\mathbb{K}^G$, so is a left unit of $\leas'(\mathbb{K} G)$.
If $e_G$ is the unit of the group $G$, then $\varepsilon=\delta_{e_G}$ is a right integral of $\mathbb{K} G$, so is a left counit of~$\leas'(\mathbb{K} G)$. As $\varepsilon(a)=1$, $\leas'(\mathbb{K} G)=\leas(\mathbb{K} G,m_\varepsilon,\Delta_a)$. For any $g,h\in G$,
\[m_\varepsilon(g\otimes h)=(\operatorname{Id} \otimes \delta_{e_G})\circ \Phi(g\otimes h)=h \delta_{e_G}\bigl(h^{-1}g\bigr)=\delta_{g,h} h.\]
For any $g\in G$,
\[\Delta_a(g)=\sum_{h\in G}\Phi(g\otimes h)=\sum_{h\in G}h\otimes h^{-1}g=\sum_{\substack{g_1,g_2\in G,\\ g_1g_2=g}}g_1\otimes g_2.\]
So $(\mathbb{K} G,m_\varepsilon,\Delta_a)$ is isomorphic to $\mathbb{K}^G$, via the map sending $g$ to $\delta_g$, for any $g\in G$.

By duality, $a$ is a left counit of $\leas'\bigl(\mathbb{K}^G\bigr)$ and $\varepsilon$ is a left unit of $\leas'\bigl(\mathbb{K}^G\bigr)$. For any $g,h\in G$,
\[m_a(\delta_g \otimes \delta_h)=(\operatorname{Id} \otimes a)\circ \Phi(g\otimes h)=\sum_{\substack{h_1,h_2\in G,\\ h_1h_2=h}} \delta_{h_1}\otimes \delta_{h_2^{-1}}\delta_g(a)=\delta_{hg}. \]
For any $g\in G$,
\[\Delta_\varepsilon(\delta_g)=\Phi(\delta_g \otimes \delta_{e_G})=\sum_{h\in G}\delta_h \otimes \delta_h \delta_g=\sum_{h\in G}\delta_h \otimes \delta_{g,h} \delta_h=\delta_g \otimes \delta_g.\]
So $\bigl(\mathbb{K}^G,m_a,\Delta_\varepsilon\bigr)$ is isomorphic to $\mathbb{K} G^{\rm op}$ via the map sending $\delta_g$ to $g$, for any $g\in G$. \end{proof}

\begin{Proposition} \label{prop5.16}\quad
\begin{enumerate}\itemsep=0pt
\item[$(1)$] The nonzero special vectors of eigenvalue $1$ of $\leas(\mathbb{K} G)$ and of $\leas'(\mathbb{K} G)$ are the elements
$\lambda \sum_{\alpha \in H} \alpha$,
where $\lambda$ is a nonzero scalar and $H$ is a subgroup of $G$.
\item[$(2)$] If $G$ is finite, the nonzero special vectors of eigenvalue $1$ of $\leas\bigl(\mathbb{K}^G\bigr)$ and of $\leas'\bigl(\mathbb{K}^G\bigr)$ are the elements
$\lambda \sum_{\alpha \in H} \delta_\alpha$,
where $\lambda$ is a nonzero scalar and $H$ is a subgroup of $G$.
\end{enumerate} \end{Proposition}

\begin{proof}
Any $a\in A$ can be written under the form $ a=\sum_{\alpha \in G} \lambda_\alpha \alpha$. Then
\begin{align*}
&a\text{ is a special vector of eigenvalue 1 of }\leas(\mathbb{K} G)\\
&\Longleftrightarrow \sum_{\alpha,\beta \in G} a_\alpha a_\beta \alpha \otimes \beta=\sum_{\alpha,\beta\in G}a_\alpha a_\beta \alpha \beta \otimes \alpha\\
&\Longleftrightarrow\sum_{\alpha,\beta \in G} a_\alpha a_\beta \alpha \otimes \beta=\sum_{\alpha,\beta\in G}a_\beta a_{\beta^{-1}\alpha} \alpha \beta \otimes \alpha\\
&\Longleftrightarrow\forall \alpha,\beta \in G, a_\beta(a_\alpha-a_{\beta^{-1}\alpha})=0.
\end{align*}
Let $a$ be a nonzero special vector of eigenvalue 1 of $\leas(\mathbb{K} G)$. Let us put $a_{1_G}=\lambda$ and $H=\{\alpha \in G, a_\alpha \neq 0\}$. Let $\alpha=\beta \in H$.
As $a_\beta\neq 0$, we obtain $a_\alpha=a_{1_G}=\lambda$, so $1_G\in H$ and~${\lambda \neq 0}$. For any $\beta \in H$, taking $\alpha=1_G$, we~obtain $a_{\beta^{-1}}=\lambda$, so $\beta^{-1}\in H$.
If $\alpha,\beta \in H$, we~obtain that $a_{\beta^{-1}\alpha}=a_\alpha\neq 0$, so $\beta^{-1}\alpha \in H$. Hence, $H$ is a subgroup and $ a=\lambda \sum_{\alpha \in H}\alpha$,
\begin{align*}
&a\text{ is a special vector of eigenvalue 1 of }\leas'(\mathbb{K} G)\\
&\Longleftrightarrow \sum_{\alpha,\beta \in G} a_\alpha a_\beta \alpha \otimes \beta=\sum_{\alpha,\beta\in G}a_\alpha a_\beta \beta \otimes \beta^{-1}\alpha\\
&\Longleftrightarrow\sum_{\alpha,\beta \in G} a_\alpha a_\beta \alpha \otimes \beta=\sum_{\alpha,\beta\in G}a_\beta a_{\alpha\beta} \alpha \beta \otimes \alpha\\
&\Longleftrightarrow\forall \alpha,\beta \in G, a_\alpha(a_\beta-a_{\alpha\beta})=0.
\end{align*}
Let $a$ be a nonzero special vector of eigenvalue 1 of $\leas'(\mathbb{K} G)$. Let us put $a_{1_G}=\lambda$ and $H=\{\alpha \in G, a_\alpha \neq 0\}$. Let $\alpha=\beta \in H$.
If $\alpha \in H$, for $\beta=1_G$, we obtain $a_{1_G}=a_\alpha=\lambda$, so $1_G\in H$ and $\lambda \neq 0$; for $\beta=\alpha^{-1}$, we obtain $a_{\alpha^{-1}}=a_{1_G}=\lambda\neq 0$, so $\alpha^{-1}\in G$.
If $\alpha,\beta \in H$, we obtain that $a_{\alpha\beta}=a_\beta\neq 0$, so $\alpha \beta\in H$. Hence, $H$ is a subgroup and $ a=\lambda \sum_{\alpha \in H}\alpha$.

Let $f\in \mathbb{K}^G$. We put $f(\alpha)=a_\alpha$ for any $\alpha \in G$,
\begin{align*}
&f\text{ is a special vector of eigenvalue 1 of }\leas\bigl(\mathbb{K}^G\bigr)\\
&\Longleftrightarrow\forall \alpha,\beta \in G, a_\alpha a_\beta=a_{\alpha\beta}\\
&\Longleftrightarrow\forall \alpha,\beta \in G, a_\alpha(a_\beta-a_{\alpha\beta})=0;\\
&f\text{ is a special vector of eigenvalue 1 of }\leas'\bigl(\mathbb{K}^G\bigr)\\
&\Longleftrightarrow\forall \alpha,\beta \in G, a_\alpha a_\beta=a_\beta a_{\beta^{-1}\alpha}\\
&\Longleftrightarrow\forall \alpha,\beta \in G, a_\beta(a_\alpha-a_{\beta^{-1}\alpha})=0.
\end{align*}
The conclusion is the same as for $\mathbb{K} G$. \end{proof}

\begin{Remark}\quad
\begin{enumerate}\itemsep=0pt
\item[(1)] From Proposition \ref{prop5.4}, the left units of $\leas(\mathbb{K} G)$ are the multiples of $e_G$, and its left counits are the multiples of its counit.
If $G$ is finite, the left units of $\leas\bigl(\mathbb{K}^G\bigr)$ are the multiple of $ \sum_{g\in G} g$,and its left counits are the multiples of $e_G$.
\item[(2)] From Proposition \ref{prop5.11}, it is not difficult to show that if $G$ is finite, the left units of $\leas'(\mathbb{K} G)$ are the multiples of $ \sum_{g\in G} g$; if $G$ is not finite, $\leas'(\mathbb{K} G)$ has no nonzero left unit.
The left counits of $\leas'(\mathbb{K} G)$ are the multiples of $\delta_{e_G}$.
By duality, if $G$ is finite, the left units of $\leas'\bigl(\mathbb{K}^G\bigr)$ are the multiples of $\delta_{e_G}$ and its left counits are the multiples of $ \sum_{g\in G} g$.
\end{enumerate}\end{Remark}

\subsection*{Acknowledgements}

The author thanks the anonymous referees for their useful and constructive comments on the first version of this text. The author acknowledges support from the grant ANR-20-CE40-0007 \emph{Combinatoire Alg\'e\-brique, R\'esurgence, Probabilit\'es Libres et Op\'erades}.

\pdfbookmark[1]{References}{ref}
\LastPageEnding

\end{document}